\documentclass[11pt,leqno,fleqn]{article}
\usepackage{amsmath,amssymb,hyperref}
\newcounter{sectie}
\newcounter{subsectie}
\newcommand{\sect}[2]{\refstepcounter{sectie}\setcounter{subsectie}{0}
\section*{\boldmath \thesectie. #2}%
\label{#1}}
\newcommand{\di}[2]{%
\refstepcounter{equation}%
\label{#1}%
\begin{list}{}{
\topsep 5mm
\leftmargin 10mm
\rightmargin 0cm
\itemsep 0mm
\listparindent 0mm
\parsep 0mm
\labelsep 1mm
\labelwidth 10mm
}%
\item[\rm (\theequation)\hfill]
\begin{list}{}{
\topsep 0mm
\leftmargin 8mm
\rightmargin 0mm
\itemsep 0mm
\listparindent 0mm
\parsep 0mm
\labelsep 1.5mm
\labelwidth 6.5mm
}
#2
\end{list}%
\end{list}%
}
\newcommand{\items}[1]{\item[#1]\vspace{-\itemsep}}
\newcommand{\rf}[1]{{\rm (\ref{#1})}}
\newcommand{\diz}[1]{%
\refstepcounter{equation}%
\begin{list}{}{
\topsep 5mm
\leftmargin 10mm
\rightmargin 0cm
\itemsep 0mm
\listparindent 0mm
\parsep 0mm
\labelsep 1mm
\labelwidth 10mm
}%
\item[\rm (\theequation)\hfill]
\begin{list}{}{
\topsep 0mm
\leftmargin 8mm
\rightmargin 0mm
\itemsep 0mm
\listparindent 0mm
\parsep 0mm
\labelsep 1.5mm
\labelwidth 6.5mm
}
#1
\end{list}%
\end{list}%
}
\newcommand{\nr}[1]{\item[{\rm (#1)}]}
\newcommand{\nrs}[1]{\item[{\rm (#1)}]\vspace{-\itemsep}}
\newcommand{\dy}[2]{%
\refstepcounter{equation}%
\label{#1}%
\begin{list}{}{
\topsep 3mm
\leftmargin 18mm
\rightmargin 0cm
\itemsep 0mm
\listparindent 0mm
\parsep 0mm
\itemsep 0mm
\labelsep 0mm
\labelwidth 18mm
}%
\item[\rm (\theequation)\hfill]
#2
\end{list}%
}
\newcommand{\oZ}{{\mathbb{Z}}}
\newcommand{\subsectz}[1]{\refstepcounter{subsectie}
\subsection*{\boldmath \thesectie.\thesubsectie. #1}%
}
\newcommand{\undef}{\text{$*$}}
\newcommand{\subsect}[2]{\refstepcounter{subsectie}
\subsection*{\boldmath \thesectie.\thesubsectie. #2}%
\label{#1}}
\newcounter{bewering}
\newcommand{\prop}[2]{\refstepcounter{bewering}\vspace{4mm}\noindent{\bf Proposition \thebewering.}\label{#1}{\it #2}}
\newcommand{\pf}{\vspace{3mm}\noindent{\bf Proof.}\ }
\newcommand{\dyyz}[1]{\dyz{\raggedright$\dps#1$}}
\newcommand{\dyz}[1]{%
\refstepcounter{equation}%
\begin{list}{}{
\topsep 3mm
\leftmargin 18mm
\rightmargin 0cm
\itemsep 0mm
\listparindent 0mm
\parsep 0mm
\itemsep 0mm
\labelsep 0mm
\labelwidth 18mm
}%
\item[\rm (\theequation)\hfill]
#1
\end{list}%
}
\newcommand{\dps}{\displaystyle}
\newcommand{\bx}{\hspace*{\fill} \hbox{\hskip 1pt \vrule width 4pt height 8pt depth 1.5pt \hskip 1pt}

\addvspace{4mm}}
\newcommand{\dist}{\text{\rm dist}}
\newcommand{\dyy}[2]{\dy{#1}{\raggedright$\dps#2$}}
\newcommand{\ac}[1]{|#1|_{\alpha}}
\newcommand{\kfrac}[2]{\mbox{$\frac{#1}{#2}$}}
\newcommand{\de}[2]{\dy{#1}{\raggedright$\displaystyle #2 $}}
\newcommand{\dez}[1]{\dyz{\raggedright$\displaystyle #1 $}}
\newcommand{\DD}{{\cal D}}
\newcommand{\VV}{{\cal V}}
\newcommand{\PP}{{\cal P}}
\newcounter{stelling}
\newcommand{\thm}[2]{\refstepcounter{stelling}\vspace{4mm}\noindent{\bf Theorem \thestelling.}\label{#1}{\it #2}}
\renewcommand{\SS}{{\cal S}}
\newcommand{\clnn}[1]{\vspace{4mm}\noindent{\em Claim.}  {\it #1}}
\newcommand{\FF}{{\cal F}}
\newcommand{\WW}{{\cal W}}
\newcommand{\pfcl}{\vspace{3mm}\noindent{\em Proof.}\ }
\newcommand{\openbx}{\hspace*{\fill} $\Box$\\ \vspace{1mm}}
\newcommand{\thmz}[1]{\refstepcounter{stelling}\vspace{4mm}\noindent{\bf Theorem \thestelling.}{\it #1}}
\newcommand{\sectz}[1]{\refstepcounter{sectie}\setcounter{subsectie}{0}
\section*{\boldmath \thesectie. #1}%
}
\renewcommand{\phi}{\varphi}
\begin{document}
\begin{center}
\baselineskip 7mm

{\large\bf\boldmath FINDING $k$ PARTIALLY DISJOINT PATHS IN A DIRECTED PLANAR GRAPH

}
\vspace{3mm}
Alexander Schrijver \footnote{ University of Amsterdam and CWI, Amsterdam.
Mailing address:
Korteweg-de Vries Institute for Mathematics, University of Amsterdam,
P.O. Box 94248, 1090 GE Amsterdam, The Netherlands.
Email: lex@cwi.nl.
[\today] ---
The research leading to these results has received funding from the European Research Council
under the European Union's Seventh Framework Programme (FP7/2007-2013) / ERC grant agreement
n$\mbox{}^{\circ}$ 339109.}

\end{center}

\noindent
{\small{\bf Abstract.}
The {\em partially disjoint paths problem} is:
{\em given:}
a directed graph, vertices $r_1,s_1,\ldots,$\linebreak[2]$r_k,s_k$, and a set $F$ of pairs
$\{i,j\}$ from $\{1,\ldots,k\}$, 
{\em find:} for each $i=1,\ldots,k$ a directed $r_i-s_i$ path $P_i$ such that if $\{i,j\}\in F$
then $P_i$ and $P_j$ are disjoint.

We show that for fixed $k$, this problem is solvable in polynomial time if the directed graph is planar.
More generally, the problem is solvable in polynomial time for directed graphs embedded on a
fixed compact surface.
Moreover, one may specify for each edge a subset of $\{1,\ldots,k\}$ prescribing which of the
$r_i-s_i$ paths are allowed to traverse this edge.

}

\sect{s1}{Introduction}

In this paper we show that the following problem,
the $k$ {\em partially disjoint paths problem}, is solvable in polynomial time for
directed planar graphs, for each fixed $k$:
\di{d1}{
\item[{\em given:}]
a directed graph $D=(V,E)$, vertices $r_1,s_1,\ldots,r_k,s_k$ of $D$, and a set $F$ of pairs
$\{i,j\}$ from $\{1,\ldots,k\}$, 
\items{\em find:} for each $i=1,\ldots,k$, a directed $r_i-s_i$ path $P_i$ in $D$ such that
if $\{i,j\}\in F$ then $P_i$ and $P_j$ are disjoint.
}
Here `disjoint' means vertex-disjoint.
So $F$ prescribes the set of pairs of paths that are forbidden to intersect.

This paper extends [15], where {\em all} pairs of paths are prescribed
to be disjoint (so $F$ is the set of all pairs from $\{1,\ldots,k\}$).
Also the method of [15] based on free groups and cohomology is extended to
free partially commutative groups (but also some simplifications of the method in
[15] have been included in the present paper).

The partially disjoint paths problem comes up in multi-commodity routing where certain commodities
are forbidden to use the same facility, to avoid clashes of conflicting commodities (radio frequencies,
soccer fan gangs, chemicals (including gases) through a pipeline network, or time slots in routing on a VLSI-chip).

The disjoint paths problem is well-studied, and generally NP-complete, implying {\em a fortiori}
that the partially disjoint paths problem is generally NP-complete.
The problem is NP-complete if we do not fix $k$, even in the undirected case (Lynch [10]).
Moreover, it is NP-complete for $k=2$ for directed graphs (Fortune, Hopcroft, and Wyllie [8]).
This is in contrast to the undirected case (if NP$\neq$P), where
Robertson and Seymour [14] showed that, for any fixed $k$, the $k$ disjoint paths 
problem is polynomial-time solvable for any graph (not necessarily planar).
For a survey of results till 2003 we refer to Chapter 70 of [16].

Our method for the partially disjoint paths problem \rf{d1} for directed planar graphs
consists of a number of layers and reductions:
\diz{
\nr{i} The top layer is to select a homology type for the solution.
The number of potentially feasible homology types can be bounded by $(2|E(G)|+1)^{4k^2}$.
This is the only level where the `fixed $k$' comes in.
\nrs{ii} For each homology type, one can find in polynomial time a solution of that type,
if it exists.
The formalism to keep track of homology is that of flows over a `graph group':
the group given by generators $g_1,\ldots,g_k$ and relations $g_ig_j=g_jg_i$ for all $i,j$ with
$\{i,j\}\not\in F$.
\nrs{iii} Finding such a solution of the prescribed homology type is done by reduction to a
`cohomology feasibility problem'
in a (generally nonplanar) {\em extension} of the planar dual of the input graph.
(This is why we need cohomology --- homology in the original, planar graph seems not enough,
mainly because disjoint paths should not only be edge-disjoint, but also vertex-disjoint.)
\nrs{iv} This cohomology feasibility problem is reduced to a 2-satisfiability problem, whose input is based
on a (polynomial) number of `pre-feasible' solutions for the cohomology feasibility problem.
\nrs{v} Finding these pre-feasible solutions forms the bottom layer of the algorithm.
It consists of a rather brute-force, but yet polynomial-time, constraint satisfaction method
(adapting an instance as long as it is not pre-feasible).
}
In our description, we start at the bottom and work our way up to the top layer.

The method rests on quite basic combinatorial group theory.
The approach allows application of the algorithm where the embedding
of the graph in the plane is given in an implicit way, viz.\ by a list of the cycles that bound
the faces of the graph, or alternatively by the clockwise order of edges incident with $v$,
for each vertex $v$.

Our method directly extends to directed graphs on
any fixed compact surface and to inputs where for each edge $e$ a subset $K_e$ of $\{1,\ldots,k\}$
is given that prescribes which of the $r_i-s_i$ paths may traverse $e$.
We did not see if the method would extend to a polynomial-time algorithm if, instead of fixing $k$,
we fix the number of faces by which $r_1,s_1,\ldots,r_k,s_k$ can be covered.

Our algorithm is a `brute force' polynomial-time algorithm.
We did not aim at obtaining the best possible running time bound, as
we presume that there are much faster (but possibly more complicated)
methods for problem \rf{d1} for directed planar graphs than the one we describe in this paper.

We could not avoid that $k$ pops up in the degree of the polynomial.
In fact, Cygan, Marx, Pilipczuk, and Pilipczuk [3] recently showed that there exists a constant $t$,
independent of $k$, such that the $k$ (fully) disjoint paths problem for directed planar
graphs is solvable in $O(n^t)$ time, for any fixed $k$.
So $k$ only shows up in the coefficient of the polynomial.
In other words, the problem is `fixed parameter tractable'.
This raises the question if also the partially disjoint paths problem is fixed parameter tractable
for directed planar graphs.

In the case of undirected planar graphs, it was shown by
Reed, Robertson, Schrijver, and Seymour [13]
that the $k$ disjoint paths problem can be solved in {\em linear} time, for any fixed $k$.
This algorithm utilizes methods from Robertson and Seymour's
graph minors theory.
For general undirected graphs, the $k$ disjoint paths problem is solvable in time $O(n^2)$ for any fixed $k$
([14], [9]).

\sect{13me03a}{Graph groups}

Our method uses the framework of combinatorial group theory,
viz.\ groups defined by generators and relations.
For background literature on combinatorial group theory we refer to 
Magnus, Karrass, and Solitar [12] and Lyndon and Schupp [11].

In particular we utilize `graph groups'.
These groups are studied {\em inter alia} by Baudisch [2],
Droms [6], Servatius [17], Wrathall [18], and Esyp, Kazachkov, and Remeslennikov [7].
Specific properties of graph groups that we will use are given in [2]
and [7], but we will
also need several other properties that seem not to have been considered before,
in particular concerning phenomena like convexity and periodicity emanating in graph groups.

We first give some standard terminology.
Let $g_1,\ldots,g_k$ form an abstract set of {\em generators}.
Call $g_1,g_1^{-1},\ldots,g_k,g_k^{-1}$ {\em symbols}.
A {\em word} (of {\em size} $t$) is a sequence $\alpha_1\cdots \alpha_t$ where each $\alpha_j$ is a symbol.
The empty word (of size 0) is denoted by $\emptyset$.
Define $(g_i^{-1})^{-1}:=g_i$, and $(\alpha_1\cdots \alpha_t)^{-1}:=\alpha_t^{-1}\cdots \alpha_1^{-1}$.

Let $g_1,\ldots,g_k$ be generators, and let $F$ be a set of unordered pairs
$\{i,j\}$ from $[k]$ with $i\neq j$.
So $([k],F)$ is an undirected graph.
(Throughout this paper: $[k]:=\{1,\ldots,k\}$.)

Then the group $G=G_F$ is generated by the generators $g_1,\ldots,g_k$, with relations%
\dy{d5a}{
$g_ig_j=g_jg_i$ for each pair $\{i,j\}\not\in F$.
}
If $F=\emptyset$p $G_F$ is the {\em free group} generated by $g_1,\ldots,g_k$.
If $F$ consists of {\em all} pairs, the group $G_F$ is the {\em free group} generated by $g_1,\ldots,g_k$.
If $F=\emptyset$ then $G_F$ is isomorphic to $\oZ^k$.
Let $1$ denote the unit element of $G_F$.
So $1=\emptyset$.

The group $G_F$ is called a {\em graph group},
or a {\em free partially commutative group},
or a {\em right-angled Artin group},
or a {\em semifree group}.
(Our definition \rf{d5a} of graph group differs in a nonessential way from that generally used, where the graph describes the pairs of commuting
generators, rather than the pairs of noncommuting generators.
Definition \rf{d5a} is more convenient for our purposes.
For instance, it implies that the group $G_F$ is equal to the product of the groups obtained from
each component of the graph $([k],F)$.)

\subsectz{Independent symbols, commuting, reduced words}

We review the basics of graph groups, referring to Baudisch [2] and Esyp, Kazachkov, and Remeslennikov [7] for the
elaboration of some details.

To describe $G_F$, call symbols $\alpha$ and $\beta$ {\em independent}
if $\alpha\in\{g_i,g_i^{-1}\}$ and $\beta\in\{g_j,g_j^{-1}\}$ for some $\{i,j\}\not\in F$ with $i\neq j$.
So if $\alpha$ and $\beta$ are independent then $\alpha\beta=\beta\alpha$ and $\beta\neq \alpha^{\pm 1}$.
(It follows from \rf{pN5b} below that also the converse implication holds.)

Call words $w$ and $v$ {\em equivalent} if $v$ if $v$ arises from $w$ by iteratively:
\di{dN5c}{
\nr{i} replacing $x\alpha\alpha^{-1}y$ by $xy$ or vice versa, where $\alpha$ is a symbol,
\nrs{ii} replacing $x\alpha\beta y$ by $x\beta\alpha y$ where $\alpha$ and $\beta$ are independent symbols.
}
By {\em commuting} we will mean applying (ii) iteratively.

Then the elements of $G_F$ are equivalence classes of words, which we can indicate by words,
although one should obviously keep in mind that different words will indicate one group element.
We will write $w\equiv v$ if we want to stress that $w$ and $v$ are equal as words.
We denote $G_F$ by $G$ if $F$ is clear from the context.

A word $w$ is called {\em reduced} if it is not equal (as a word) to
$x\alpha y\alpha^{-1}z$ for some symbol $\alpha$ independent of $y$.
Note that reducedness is a property of words, and that it is invariant under commuting.
We say that a symbol $\alpha$ {\em occurs in} an element $x$ of $G$, or that $x$ {\em contains} $\alpha$,
if $\alpha$ occurs in any reduced word representing $x$.
Two elements $x$ and $y$ of $G$ are called {\em independent} if any symbol in $x$ and any symbol
in $y$ are independent.
(In particular, $\beta\neq \alpha^{\pm 1}$ for any symbols $\alpha$ in $x$ and $\beta$ in $y$.)

The following is basic --- see Lemma 2.3 in [7]:
\dy{pN5b}{\em
Let $w$ and $x$ be reduced words with $w=x$ as group elements.
Then word $x$ can be obtained from $w$ by a series of commutings.
}
Define, for $x\in G$, $|x|$ as the size of any reduced word representing $x$.
So $|xy|\leq|x|+|y|$ for all $x,y\in G$.

\rf{pN5b} implies that testing if $w=1$ is easy:
just replace (iteratively) any contiguous subword
$\alpha y\alpha^{-1}$ by $y$ where $\alpha$ is a symbol and $y$ is a word independent of $\alpha$.
The final word is empty if and only if $w=1$.
This gives a test for equivalence of words $w$ and $x$: just test if $wx^{-1}=1$. 
So the `word problem' for graph groups is solvable in polynomial time.
(In fact it can be solved in linear time --- see Wrathall [18].)

It will be convenient to have a mean to emphasize when the concatenation of two reduced words $x$ and $y$ gives again a
reduced word (without cancellation as in \rf{dN5c}(i)).
In other words, when $|xy|=|x|+|y|$.

To this end, we add an abstract new element $\undef$ to $G$ and define a multiplication $\cdot$ on $G\cup\{\undef\}$ as follows.
Let $x,y\in G$.
Then $x\cdot y:=xy$ if $|xy|=|x|+|y|$,
and $x\cdot y:=\undef$ if $|xy|<|x|+|y|$.
So $x\cdot y$ belongs to $G$ if for any reduced words $x'$ and $y'$ representing $x$ and $y$ one has
that the concatenation of $x'$ and $y'$ is reduced.
So no symbol in $x'$ cancels out any symbol in $y'$.
If we moreover set $\undef\cdot x:=\undef$ and $x\cdot\undef:=\undef$ for all $x\in G\cup\{\undef\}$,
we obtain an associative multiplication $\cdot$ on $G\cup\{\undef\}$.

The only purpose of introducing $*$ is to have a convenient and formally correct tool to write,
for $x,y,z\in G$, $x=y\cdot z$, meaning $x=yz$ and $|x|=|y|+|z|$.
By extension, for $x,y_1,\ldots,y_n\in G$,
$x=y_1\cdot\ldots\cdot y_n$ means $x=y_1y_2\ldots y_n$ and $|x|=|y_1|+|y_2|+\cdots+|y_n|$.
That is, in the concatenation of reduced words $y_1,\ldots,y_n$ there is no cancelation.
The element $\undef$ will not occur anymore below.

While the multiplication $\cdot$ is associative,
it is generally not the case that if $xy=x\cdot y$ and $yz=y\cdot z$ then $xyz=x\cdot y\cdot z$,
because in $xyz$, symbols in $x$ might cancel out symbols in $z$.
Nevertheless, the following holds.
Call $y\in G$ a {\em segment} of $a\in G$ if there exist $x,z\in G$ with $a=x\cdot y\cdot z$.
Then:
\dy{29no14f}{\em
If $xy=x\cdot y$ and $yz=y\cdot z$, then $y$ is a segment of $xyz$.
}
To see this, let $x$, $y$, and $z$ be reduced words, and consider the concatenation of $x,y,z$.
In the cancellation to obtain a reduced word, only symbols in $x$ and symbols in $z$ can cancel
each other out (since the concatenations $x,y$ and $y,z$ are reduced).
So $y$ survives as segment of $xyz$.

\subsect{s12A}{The partial order $\leq$}

Most of what follows in this section is, explicitly or implicitly, in [7].
Let $x,y\in G$.
We write $x\leq y$ if $x=x\cdot a$ for some $a$ (namely, $a=y^{-1}x$).
(If $x\leq y$, $x$ is called in [7] a {\em left divisor} of $y$.)
So if $y$ is given as reduced word, it means that $y$ can be commuted so that the first
$|x|$ symbols in $y$ form $x$ (by \rf{pN5b}).
Also, $x\leq y$ if and only if $|y|=|x|+|x^{-1}y|$.

It is easy to derive from the norm properties of $|.|$ that $\leq$ is a
partial order.
In fact, the partial order $\leq$ is a lattice if we add to $G$ an element
$\infty$ at infinity (Propositions 3.10 and 3.12 in [7]).
This follows from the existence of the meet $x\wedge y$ for all $x,y\in G$.
Then the join $\vee$ exists for all $x,y$ for which there exists $z\in G$ with $x,y\leq z$
(then $x\vee y$ is the meet of all such $z$).
So adding an element $\infty$ with $\infty\geq x$ for all $x$,
makes $(G\cup\{\infty\},\leq)$ to a lattice.
Then $x\vee y=\infty$ if there is no $z\in G$ with $z\geq x,y$.

If finite, the join $x\vee y$ can be described as follows (Proposition 3.18 in [7]):
\dy{1no14b}{\em
Let $x,y\in G$ and define $x':=(x\wedge y)^{-1}x$ and $y':=(x\wedge y)^{-1}y$.
Then $x\vee y<\infty$ if and only if $x'$ and $y'$ are independent.
}
Moreover:
\dy{1de14b}{\em
If $x\vee y<\infty$ then $x\vee y=(x\wedge y)\cdot x'\cdot y'=x(x\wedge y)^{-1}y$.
}

For any $a\in G$ define:
\dy{d27g}{
$a^{\downarrow}:=\{x\in G\mid x\leq a\}$ and $a^{\uparrow}:=\{x\in G\mid x\geq a\}$.
}
The norm characterization of $\leq$ implies for all $x,y,z$ with $x\leq y,z$:
\dy{d26bi}{\em
$y\leq z$ if and only if $x^{-1}y\leq x^{-1}z$.
}
Hence for any $a\in G$,
the function $a^{\uparrow}\to a^{-1}a^{\uparrow}$ with $b\mapsto a^{-1}b$ for $b\in a^{\uparrow}$
is an order isomorphism, and therefore:
\dy{25ja15d}{\em
If $b,c\geq a$ then $a^{-1}(b\wedge c)= a^{-1}b\wedge a^{-1}c$, and if moreover $b\vee c<\infty$,
then $a^{-1}(b\vee c)=a^{-1}b\vee a^{-1}c$.
}

\prop{p5c}{
Let $x_1,\ldots,x_t$ be such that $x_i\vee x_j<\infty$ for all $i,j$.
Then $x_1\vee\cdots\vee x_t<\infty$.
}

\pf
It suffices to show this for $t=3$.
Suppose $x\vee y\vee z=\infty$ while $y\vee z<\infty$.
Let $a:=x\wedge(y\vee z)$.
As $x\vee y\vee z=\infty$, $a^{-1}x$ and $a^{-1}(y\vee z)$ are dependent.
Now, using
\rf{25ja15d} (with $b:=a\vee y$ and $c:=a\vee z$) and \rf{1de14b}:
\dyyz{
a^{-1}(y\vee z)=
a^{-1}(a\vee y\vee z)=
a^{-1}(a\vee y)\vee a^{-1}(a\vee z) =
(a\wedge y)^{-1}y\vee(a\wedge z)^{-1}z=
(x\wedge y)^{-1}y\vee(x\wedge z)^{-1}z.
}
So $a^{-1}x$ is dependent on at least one of $(x\wedge y)^{-1}y$ and $(x\wedge z)^{-1}z$,
say $a^{-1}x$ is dependent on $(x\wedge y)^{-1}y$.
Since $x\wedge y\leq a$ this implies that $(x\wedge y)^{-1}x$ is dependent $(x\wedge y)^{-1}y$.
So $x\vee y$ is infinite.
\bx

We say that a symbol $\alpha$ is a {\em minimal symbol} of $x\in G$ if $\alpha\leq x$
(so $\alpha$ can be commuted so as to become the first symbol).
Similarly, $\alpha$ is a {\em maximal symbol} if $\alpha^{-1}\leq x^{-1}$, i.e.,\ if
$x\alpha^{-1}\leq x$
(so $\alpha$ can be commuted so as to become the last symbol).

We will need that for all $x,y,z\in G$:
\dy{d26bii}{\em
If $y\leq xyz=x\cdot y\cdot z$, then $y\leq xy$.
}
This can be seen by induction on $|z|$.
Let $\alpha$ be a maximal element of $z$, let $z':=z\alpha^{-1}$, and suppose
that $y\not\leq xyz'$.
Let $\alpha$ occur $m$ times in $y$.
As $y\leq xyz$ and $y\not\leq xyz'$, $\alpha$ occurs $m-1$ times in $xyz'$.
This contradicts the fact that $\alpha$ occurs $m$ times in $y$,
hence in $x\cdot y\cdot z'$.
This proves \rf{d26bii}.

We also will need that for all $x,y,z\in G$:
\dy{25ja15b}{\em
If $x,y\leq z$ then $(x\wedge y)^{-1}x\leq y^{-1}z$.
}
To see this, let $b:=x\vee y=y(x\wedge y)^{-1}x$.
As $y\leq b\leq z$, we have $(x\wedge y)^{-1}x=y^{-1}b\leq y^{-1}z$
by \rf{d26bi}.
proving \rf{25ja15b}.

Moreover, for all $x,y,z\in G$:
\dy{31ok14b}{\em
If $y^{-1}x\wedge y^{-1}z=1$ then $(x\wedge z)\vee(x\wedge y)\vee(y\wedge z)=y$.
}
To prove this, we can assume (by \rf{25ja15d}) that $x\wedge y\wedge z=1$.
Let $a:=y\wedge z$, $b:=x\wedge z$, and $c:=x\wedge y$.
Then $a\wedge b=a\wedge c=b\wedge c=1$.
Hence, by \rf{1no14b},
$a$, $b$,and $c$ are pairwise independent, and $bc\leq x$, $ac\leq y$, $ab\leq z$.
Let $x'$, $y'$, and $z'$ satisfy $x=bcx'$, $y=acy'$, and $z=abz'$.
Since $x\wedge y=c$, we know $y^{-1}x={y'}^{-1}\cdot a^{-1}\cdot b\cdot x'$.
Hence, as $a$ and $b$ are independent, ${y'}^{-1}\cdot b\leq y^{-1}x$.
Similarly, ${y'}^{-1}\cdot b\leq y^{-1}z$.
As $y^{-1}x\wedge y^{-1}z=1$, we have $y'=1$ and $b=1$.
Hence $(x\wedge z)\vee(x\wedge y)\vee(y\wedge z)=a\vee b\vee c=ac=y$.

(We finally remark, but will not need, that
the lattice $G\cup\{\infty\}$ is not distributive (if $F\neq\emptyset$), while for each $a\in G$,
the sublattice $a^{\downarrow}$ is distributive.)

\subsect{s26a}{Convex sets}

The function $\dist(x,y):=|x^{-1}y|$ is a metric, since, for all $x,y\in G$,
(i) $|x|=0\iff x=1$, (ii) $|x^{-1}|=|x|$, and (iii) $|xy|\leq |x|+|y|$.
Note that this distance is left-invariant: $\dist(zx,zy)=\dist(x,y)$ for all $x,y,z$.

We call a subset $L$ of $G$ {\em convex} if $L$ is nonempty and
if $x,z\in L$ and $\dist(x,y)+\dist(y,z)=\dist(x,z)$ then $y\in L$.
Since the distance function is left-invariant,
if $L$ is convex also $yL$ is convex, for each $y\in G$.

\prop{p3a}{
A nonempty subset $L$ of $G$ is convex if and only if
\di{d3aa}{
\nr{i} if $x\leq y\leq z$ and $x,z\in L$ then $y\in L$,
\nrs{ii} if $x,y\in L$ then $x\wedge y\in L$ and, if $x\vee y$ is finite,
$x\vee y\in L$.
}
}

\pf
Necessity follows from the facts that if $x\leq y\leq z$ then $\dist(x,y)+\dist(y,z)=\dist(x,z)$,
that $\dist(x,y)=\dist(x,x\wedge y)+\dist(x\wedge y,y)$ and that,
if $x\vee y$ is finite then $\dist(x,y)=\dist(x,x\vee y)+\dist(x\vee y,y)$.

To see sufficiency, let $\dist(x,y)+\dist(y,z)=\dist(x,z)$ with $x,z\in L$.
We must show $y\in L$.
So $|x^{-1}y|+|y^{-1}z|=|x^{-1}z|$, hence $y^{-1}x\wedge y^{-1}z=1$.
Hence by \rf{31ok14b}, $(x\wedge z)\vee(x\wedge y)\vee(y\wedge z)=y$.
So $x\wedge z\leq x\wedge y\leq x$ and $x\wedge z\leq y\wedge z\leq z$.
Therefore, $x\wedge y$ and $y\wedge z$ belong to $L$ and hence $y$ belongs to $L$.
\bx

This implies that $a^{\uparrow}$ and $a^{\downarrow}$ are convex.
Moreover, each convex set $L$ has a unique minimal element $\min L$.
This in fact characterizes convex sets:
\dy{25ok14g}{\em
A nonempty subset $L$ of $G$ is convex if and only if for each $a\in G$,
$aL$ has a unique minimal element.
}
Here necessity follows from \rf{d3aa}(ii).
To see sufficiency, let $x,z\in L$ and $y\in G$ such that $\dist(x,y)+\dist(y,z)=\dist(x,z)$.
We prove $y\in L$.
We may assume $y=1$ (as the condition is invariant under resetting $L\to y^{-1}L$).
Let $a$ be the unique minimal element in $L$.
So $x\geq a$ and $z\geq a$.
On the other hand, $|x|+|z|=|x^{-1}z|$, and hence $x\wedge z=1$.
So $a=1$, proving \rf{25ok14g}.

Clearly, the intersection of any number of convex sets is convex again.
Moreover, convex sets satisfy the following `Helly-property':
\dy{pj4a}{\em
Let $L_1,L_2,L_3$ be convex sets with $L_i\cap L_j\neq\emptyset$ for all $i,j=1,2,3$.
Then $L_1\cap L_2\cap L_3\neq\emptyset$.
}
For choose $x\in L_1\cap L_2, y\in L_1\cap L_3, z\in L_2\cap L_3$.
Without loss of generality, $z=1$ (as we can replace $L_1,L_2,L_3$ by $z^{-1}L_1,z^{-1}L_2,z^{-1}L_3$).
Now $x\wedge y\in L_1\cap L_2\cap L_3$.

This proves \rf{pj4a}, which implies the following.
As usual, define $XY:=\{xy\mid x\in X,y\in Y\}$
and $X^{-1}:=\{x^{-1}\mid x\in X\}$, for $X,Y\subseteq G$.
Then:
\dy{p1a}{\em
Let $L_1$, $L_2$, and $L_3$ be convex, with $L_1\cap L_2\neq\emptyset$.
Then $L_1L_3^{-1}\cap L_2L_3^{-1}=(L_1\cap L_2)L_3^{-1}$.
}
Indeed, trivially, $L_1L_3^{-1}\cap L_2L_3^{-1}\supseteq(L_1\cap L_2)L_3^{-1}$.
To see the reverse inclusion, let $x\in L_1L_3^{-1}\cap L_2L_3^{-1}$.
Since $x\in L_1L_3^{-1}$, we know $x^{-1}L_1\cap L_3\neq\emptyset$.
Similarly, $x^{-1}L_2\cap L_3\neq\emptyset$.
Since also $L_1\cap L_2\neq\emptyset$, \rf{pj4a} gives
$x^{-1}L_1\cap x^{-1}L_2\cap L_3\neq\emptyset$.
Hence $x\in (L_1\cap L_2)L_3^{-1}$.

\subsectz{Ideals and closed sets}

A subset $I$ of $G$ an {\em ideal} if $I$ is nonempty and
\di{d3ba}{
\nr{i} if $y\leq x$ and $x\in I$ then $y\in I$,
\nrs{ii} if $x,y\in I$ and $x\vee y$ is finite then $x\vee y\in I$.
}
Since $1$ belongs to any ideal,
by comparing \rf{d3aa} and \rf{d3ba} one sees that each ideal is convex.
Moreover, for any $L\subseteq G$ and $x\in L$:
\dy{25ok14f}{\em
$L$ is convex if and only if $x^{-1}L$ is an ideal.
}

Call $H$ {\em closed} if both $H$ and $H^{-1}$ are ideals.
In particular, if $H$ is closed and $x\in H$, then any segment of $x$ belongs to $H$.

\prop{24ok14c}{
If $H$ is closed and $x,y\in G$, then
$x\leq y$ implies $\min(xH)\leq\min(yH)$.
}

\pf
We can assume that $y=x\alpha$ for some symbol $\alpha$.
Let $c\in H$ with $yc=\min(yH)$.
It suffices to show that there exists $d\in H$ with $xd\leq yc$.

If $\alpha^{-1}\leq c$, let $d:=\alpha c$, in which case $d\in H$ (as $H^{-1}$ is an ideal) and $xd=yc$.
If $\alpha^{-1}\not\leq c$, then let $d:=c$, in which case $xd\leq yc$.
Indeed, as $yc=\min(yH)$ and $H$ is an ideal, we know $c\leq y^{-1}$.
Since $c\leq y^{-1}$ and $\alpha^{-1}\leq y^{-1}$, $\alpha^{-1}$ and $c$ are independent.
As moreover $\alpha^{-1}\leq y^{-1}$, we know $\alpha^{-1}\leq c^{-1}y^{-1}$.
Therefore, $yc\alpha^{-1}\leq yc$, and hence $xc=y\alpha^{-1}c=yc\alpha^{-1}\leq yc$.
\bx

This is used in showing:
\dy{p7f}{\em
If $L$ is convex and $H$ is closed, then $LH$ is convex.
Moreover, $\min(LH)=\min(\min(L)H)$.
}
To show that $LH$ is convex,
by \rf{25ok14g} it suffices to show that $LH$ has a unique minimal element
(as each $xL$ is again convex).
Let $a=\min(L)$ and choose $c\in H$ with $ac:=\min(aH)$.
Then for each $x\in L$ and $y\in H$, we have by Proposition \ref{24ok14c},
as $a\leq x$, $ac=\min(aH)\leq\min(xH)\leq xy$.
So $ac$ is the unique minimal element in $LH$.
Hence $\min(LH)=\min(aH)=\min(\min(L)H)$, and we have \rf{p7f}.

This implies:
\dyz{\em
If $H$ and $H'$ are closed, then $HH'$ is closed.
}
Indeed, $H$ is an ideal, hence convex, hence by \rf{p7f}, $HH'$ is convex.
As $1\in HH'$, it follows that $HH'$ is an ideal.
Similarly, $(HH')^{-1}$ is an ideal.
So $HH'$ is closed.

This gives for any closed $H$ and $x,x\in G$:
\dyy{15de14a}{
xHz^{-1}=
x^{\downarrow}H(z^{\downarrow})^{-1}
\cap
x^{\uparrow}H(z^{\downarrow})^{-1}
\cap
x^{\downarrow}H(z^{\uparrow})^{-1}
\cap
x^{\uparrow}H(z^{\uparrow})^{-1}.
}

This follows from \rf{p1a} and \rf{p7f},
as $x^{\uparrow}$, $x^{\downarrow}$,
$z^{\uparrow}$, and $z^{\downarrow}$ all are convex,
hence $z^{\downarrow}H^{-1}$ and $z^{\uparrow}H^{-1}$ are convex.
Then \rf{p1a} gives that the right-hand side in \rf{15de14a} is equal to
$xH(z^{\downarrow})^{-1}\cap xH(z^{\uparrow})^{-1}$.
Applying \rf{p1a} to the inverse of this set, we obtain \rf{15de14a}.

\subsectz{ Peaks}

An element of $G$ is called a {\em peak} if it has precisely one maximal symbol.
The peaks are precisely the join-irreducible elements of $G$ with respect to $\vee$.
If the maximal symbol equals $\alpha$, then $p$ is called an {\em $\alpha$-peak}.
For each $x\in G$ and symbol $\alpha$, all $\alpha$-peaks $p\leq x$ are totally ordered by $\leq$.
Moreover:
\dy{16no14a}{\em
Each $x\in G$ is the join of all peaks $p\leq x$.
}
To see this, let $y$ be the join of all peaks $p\leq x$.
If $y\neq x$, choose a maximal symbol $\alpha$ of $y^{-1}x$.
Then $\alpha$ is also a maximal symbol of $x$.
Write $x=\xi_1\cdot\ldots\cdot\xi_n$ with symbols $\xi_1,\ldots,\xi_n$, in such a way that
the maximum $j$ for which $\xi_j=\alpha$ is minimized.
Then $q:=\xi_1\cdot\ldots\cdot\xi_j$ is an $\alpha$-peak with $q\leq x$.
So $q\leq y$.
This however contradicts the fact $y^{-1}x$ has maximal symbol $\alpha$,
thus showing \rf{16no14a}.

For each $x\in G$ and symbol $\alpha$, let $\ac{x}$ be the number of occurrences of symbol $\alpha$
in $x$ (not considering $\alpha^{-1}$).

\prop{13ma15h}{
Let $x,y\in G$ and let $p\leq x$ and $q\leq y$ be $\alpha$-peaks satisfying
$\ac{p^{-1}x}=\ac{q^{-1}y}$.
Then $|pq^{-1}|\leq |xy^{-1}|$.
}

\pf
By induction on $|x|+|y|$.
If $x$ is not an $\alpha$-peak, let $\beta$ be a maximal symbol of $x$
with $\beta\neq\alpha$.
Let $x':=x\beta^{-1}$.
If $\beta$ is also a maximal symbol of $y$, then
we can apply induction to $x'$ and $y':=y\beta^{-1}$, since
$x'(y')^{-1}=xy^{-1}$.
If $\beta$ is not a maximal symbol of $y$ then $|x'y^{-1}|\leq|xy^{-1}|$
(as $\beta$ is not canceled in the concatenation of $x$ and $y^{-1}$),
and hence we can apply induction to $x'$ and $y$.

So we can assume that $x$ is an $\alpha$-peak, and similarly that
$y$ is an $\alpha$-peak.
If $\ac{p^{-1}x}=0$, then $x=p$ and $y=q$, and we are done.
If $\ac{p^{-1}x}>0$, then
$x>p$ and $y>q$ and we can apply induction to
$x':=x\alpha^{-1}$ and $y':=y\alpha^{-1}$.
Note that $\ac{p^{-1}x'}=\ac{p^{-1}x}-1=\ac{q^{-1}y}-1=\ac{q^{-1}y'}$.
\bx

We use this proposition only in obtaining the following:
\dy{15ma15f}{\em
Let $p\leq r\leq ar$ with $p$ an $\alpha$-peak.
Then there exists $a'$ with $p\leq a'p$, $\ac{p^{-1}a'p}=\ac{r^{-1}ar}$, and $|a'|\leq |a|$.
}
This follows by applying
Proposition \ref{13ma15h} to $x:=r$ and $y:=ar$, taking for $q$ the (unique) $\alpha$-peak satisfying
$p\leq q\leq ar$ with $\ac{q^{-1}ar}=\ac{p^{-1}r}$, which shows that we can
take $a':=qp^{-1}$.
(Note that $\ac{q^{-1}ar}=\ac{p^{-1}r}$ is equivalent to
$\ac{p^{-1}q}=\ac{r^{-1}ar}$, since
$\ac{p^{-1}q}+\ac{q^{-1}ar}=\ac{p^{-1}ar}=\ac{p^{-1}r}+\ac{r^{-1}ar}$.)

We also need:
\dy{9fe15b}{\em
If $x\leq y$ and $\alpha$ is a symbol not occurring in $x^{-1}y$,
then $|xp^{-1}|\leq|yp^{-1}|$ for each $\alpha$-peak $p$.
}
Indeed, by induction we can assume that $x^{-1}y=\beta$ for some symbol $\beta\neq\alpha$.
Then in $yp^{-1}$, the maximal symbol $\beta$ of $y$ is not cancelled, since otherwise $\beta$ would be
maximal symbol also of $p$, hence $\beta=\alpha$, contradicting our assumption.
Hence $|xp^{-1}|\leq|x\beta p^{-1}|=|yp^{-1}|$.

\subsect{s12C}{Connectedness and cyclic reducedness}

We study periodicity of symbols in elements of $G$ in order to obtain control on `stalling' in the algorithm.
For this we need Proposition \ref{26no14b} below --- the other results in Sections
\ref{13me03a}.\ref{s12C}--\ref{13me03a}.\ref{27ja15b} are only needed to prove
Proposition \ref{26no14b}.

Call an element $b$ of $G$ {\em connected} if the generators occurring in $b$ induce a connected
subgraph of $([k],F)$.
So $b$ is connected if and only if there are no $a,c\in G$ with $b=ac$, $a\neq 1\neq c$, and
$a$ and $c$ independent.
Each peak is connected.

Call an element $b$ of $G$ {\em cyclically reduced} if $b\wedge b^{-1}=1$.
So $b$ is cyclically reduced if and only if $b^2=b\cdot b$.
Also, if $b$ is cyclically reduced, then for each $s\geq 0$: $b^s=b\cdot b\cdot\ldots\cdot b$
(cf.\ [2]).

The following proposition will be used in proving Propositions
\ref{1fe15c} and \ref{401b}.

\prop{3fe15c}{
Let $c,d\in G$ satisfy $d\leq dc$, $c\leq dc$, and $d\not\leq c$.
Suppose that $d$ is connected and that all minimal symbols of $c$ occur in $d$.
Then $|c|\leq|c\wedge d|^2$.
If $d$ is moreover cyclically reduced, then $c\leq d^{|c\wedge d|}$.
}

\pf
The proof is by induction on $|c|$.
Let $c':=(c\wedge d)^{-1}c$.
Then $c'\leq c$, by \rf{25ja15b}, since $c,d\leq dc$.
Also, $d\leq dc'$ (as $d\leq dc$ and $c'\leq c$)
and $c'\leq dc'$ (by \rf{d26bii}, as $c'\leq c\leq dc$ and $dc'\leq dc$).
As $d\not\leq c$ and $c'\leq c$, we know $d\not\leq c'$.
Moreover, as $c'\leq c$, all minimal symbols of $c'$ occur in $c$, hence in $d$.

If $c'\wedge d<c\wedge d$, then by induction $|c'|\leq |c'\wedge d|^2$ and, if $d$ is
cyclically reduced, $c'\leq d^{|c'\wedge d|}$.
Hence
$|c|=|c\wedge d|+|c'|\leq |c\wedge d|+|c'\wedge d|^2\leq|c\wedge d|^2$ and
$c\leq c\vee d=d(c\wedge d)^{-1}c=dc'\leq d^{|c'\wedge d|+1}\leq d^{|c\wedge d|}$,
as required.

So we can assume $c'\wedge d=c\wedge d$.
As $c,d\leq dc$, $c\vee d<\infty$.
So $(c\wedge d)^{-1}c=c'$ and $(c\wedge d)^{-1}d$ are independent.
Hence $c'\wedge d=c\wedge d$ and $(c\wedge d)^{-1}d$ are independent.
As $c\wedge d\neq d$ and as $d$ is connected, we know $c\wedge d=1$.
So $c$ and $d$ are independent.
Since all minimal symbols of $c$ occur in $d$, this implies $c=1$ and the bounds are trivial.
\bx


\subsectz{Conjugates}

An element $c$ of $G$ is called a {\em conjugate} of $a\in G$ if $c=x^{-1}ax$
for some $x\in G$.
Then:
\dy{15ma15y}{\em
For each $a\in G$, each conjugate $c$ of $a$ contains a segment
$x^{-1}ax$ with $x$ using only generators occurring in $a$.
}
Indeed, let $c=y^{-1}ay$.
Then \rf{15ma15y} can be proved by induction on $|y|$.
If $z:=y\wedge a\neq 1$, replace $y$ by $z^{-1}y$ and
$a$ by $z^{-1}az$, and apply induction (this resetting does not change $y^{-1}ay$).
So we can assume that $y^{-1}a=y^{-1}\cdot a$ and similarly that
$ay=a\cdot y$.
Hence $a$ is a segment of $y^{-1}ay$ by \rf{29no14f},
proving \rf{15ma15y}.

We use this in proving:
\dy{15ma15x}{\em
If $a$ and $b$ are independent, then each conjugate $c$ of $ab$ contains a segment
which is a conjugate of $a$.
}
Indeed, by \rf{15ma15y} $c$ has a segment $x^{-1}abx$ with
$x$ only using generators occurring in $ab$.
As $a$ and $b$ are independent we can write $x=yz$ with $y$ only using
generators occurring in $a$ and $z$ only using generators occurring in $b$.
Hence $y^{-1}ay$ and $z^{-1}bz$ are independent, and so
$y^{-1}ay$ is a segment of $x^{-1}abx$.

Proposition \ref{3fe15c} implies:

\prop{1fe15c}{
Let $d$ be connected and cyclically reduced.
Then for each $n\geq 0$, each conjugate $c$ of $d^{n+2|d|}$ contains $d^n$ as segment.
}

\pf
Choose $x\in G$ with $c=x^{-1}d^{n+2|d|}x$ and $|x|$ as small as possible.
Then $d\not\leq x$, otherwise replacing $x$ by $d^{-1}x$ contradicts the minimality of $x$.
Let $y:=x\wedge d^{n+|d|}$.
Then $y\leq d^{n+|d|+1}$ and $d\leq dy\leq d^{n+|d|+1}$, hence by
\rf{d26bii}, $y\leq dy$.
Since $y\leq d^{n+|d|}$, all minimal symbols of $y$ occur in $d$.
Hence by Proposition \ref{3fe15c}, $y\leq d^{|d|}$.
So $d^{|d|}=y\cdot a$ for some $a$.
Hence for $z:=y^{-1}x$ one has $x^{-1}d^{n+|d|}=z^{-1}\cdot a\cdot d^n$,
implying $x^{-1}d^{n+|d|}=x^{-1}d^{|d|}\cdot d^n$.

By symmetry, $d^{n+|d|}x=d^n\cdot d^{|d|}x$.
So by \rf{29no14f}, $d^n$ is a segment of
$x^{-1}d^{|d|}d^nd^{|d|}x=c$.
\bx

\subsectz{Periodicity}

We give conditions for the eventual periodicity of a peak:

\prop{401b}{
Let $q$ be connected and contain symbol $\alpha$,
and let $p$ be an $\alpha$-peak with $p\leq pq$.
Then there exists an $\alpha$-peak $r$ and $t\geq 0$ with $p=r\cdot q^t$
and $|r|\leq2|pqp^{-1}|^2$.
}

\pf
Let $a:=pqp^{-1}$.
Then $|q|\leq|a|$, as $|p|+|q|=|pq|=|ap|\leq|p|+|a|$.
If $|p|\leq 2|a|^2$, we can take $r:=p$ and $t:=0$.
So we can assume that $|p|>2|a|^2$.

Let $m:=\ac{q}$, and let
$p'$ be minimal with the properties that $p'\leq p$ and $\ac{p'}\geq\ac{p}-m$.
Then $p'=1$ or $p'$ is an $\alpha$-peak.
By showing that $(p')^{-1}p=q$ we are done, since then we can apply induction, as
$p'q(p')^{-1}=pqp^{-1}$.

Define $c:=p^{-1}$, $d:=q^{-1}$, and $u:=c\wedge dc$.
As $u\leq c$, $u\leq dc$, and $d\leq dc$, we know $u\leq du=d\cdot u$ (by
\rf{d26bii}).
Since $a=(dc)^{-1}c$, we have $|a|=|c|+|dc|-2|u|$, and so, as $|c|=|p|>2|a|^2$ and $|a|\geq|d|$:
\dyyz{
2|u|=|c|+|dc|-|a|=2|c|+|d|-|a|>4|a|^2+|d|-|a|\geq 2|d|^2+2|d|.
}
Hence $|u|>|d|^2+|d|\geq |u\wedge d|^2$.
Therefore, $d\leq u$ by Proposition \ref{3fe15c}.

Let $b:=p^{-1}p'$.
So we must show $b=d$.
Since $q^{-1}=d\leq u\leq c=p^{-1}$ and $m=\ac{q}$, we know $p'\leq pq^{-1}$,
that is, $d\leq b$.
On the other hand,
$b\wedge u\leq d$, since $b\wedge u\leq u\leq dc=q^{-1}p^{-1}$,
$\ac{b\wedge u}\leq m=\ac{q}$, and $p$ is a peak.
So $b\wedge u=d$, and hence we must show $b\leq u$.

Let $u':=d^{-1}u$.
Since $d\leq u\leq du$, we have $u'\leq u$ (by \rf{d26bi}).
Moreover, as $b,u\leq c$, $b\vee u<\infty$.
Hence $u'\leq u\leq b\vee u=bd^{-1}u=b\cdot u'$.
Since $|u'|=|u|-|d|>|d|^2=|u\wedge b|^2\geq|u'\wedge b|^2$, Proposition \ref{3fe15c}
gives $b\leq u'$, implying $b\leq u$.
\bx

A g.c.d argument shows:
\dy{28ja15ax}{\em
Let $p,r,r',q,q'$ be $\alpha$-peaks, and let $t,t'\geq 0$ with $p=r\cdot q^t=r'\cdot(q')^{t'}$
and $\kfrac13\ac{p}\geq |r|,|r'|,|q|,|q'|$.
Then there exists $d$ such that $q$ and $q'$ are powers of $d$.
}
Indeed, define $g:=\max\{|r|,|r'|,|q|,|q'|\}$,
$m:=\ac{p}$, $u:=\ac{q}$, and $u':=\ac{q'}$.
So $m\geq 3g$.
As $p$ is an $\alpha$-peak, we can uniquely write
$p=p_1\cdot p_2\cdot\ldots\cdot p_m$,
with each $p_i$ being an $\alpha$-peak.
Since $r$ is an $\alpha$-peak with $\ac{r}\leq g$,
and since $p=r\cdot q^t$,
the sequence $z:=(p_{g+1},\ldots,p_{m})$ is periodic with period $u$.
Similarly, $z$ is periodic with period $u'$.
Moreover, $z$ has at least $u+u'$ terms, since $m-g\geq 2g\geq u+u'$.
This implies\footnote{
If $a$ and $b$ are periods of $x=(x_1,\ldots,x_n)$ and $n\geq a+b$,
then $a-b$ is a period of $x$:
if $i\leq n-a$ then $x_{i+(a-b)}=x_{(i+a)-b}=x_{i+a}=x_i$;
if $n-a<i\leq n-(a-b)$ then $i>b$, hence $x_{i+(a-b)}=x_{(i-b)+a}=x_{i-b}=x_i$.
}
that $z$ is periodic with period $v:=\gcd\{u,u'\}$.
Let $d:=z_{m-v+1}\cdot\ldots\cdot z_m$.
Then $q=z_{m-u+1}\cdot\ldots\cdot z_m=d^{u/v}$ and similarly $q'=d^{u'/v}$.

\subsect{27ja15b}{A main tool}

We now come to a main tool for bounding the complexity of our algorithm
(which we will use in Section \ref{13me03b}.\ref{s5a}).

\prop{26no14b}{
Let $p$ be an $\alpha$-peak, and let $a,a'\in G$ be such that $p\leq ap$ and $p\leq a'p$ and
such that $\alpha$ occurs in ${p}^{-1}a'p$.
If $|p|\geq 8|a|^3$, then each conjugate of $a$ has a segment $s$ satisfying
$\ac{p^{-1}aps^{-1}}\leq 2|a'|^2$.
}

\pf
We can assume $\ac{p^{-1}ap}>2|a'|^2$, as otherwise we can take $s:=1$.
This implies $|a|\geq |ap|-|p|=|p^{-1}ap|\geq|a'|$.

Let $q$ be the component of $p^{-1}ap$ that contains $\alpha$;
that is, $q$ is the element such that $p^{-1}ap=qu$ for some $u$ independent of $q$, with $q$
connected and containing $\alpha$.
Similarly, let $q'$ be the component of $p^{-1}a'p$ that contains $\alpha$.
Note that $|q|\leq|a|$ and $|q'|\leq|a'|$.

By \rf{9fe15b} applied to $x:=pq$ and $y:=ap$ we have $|pqp^{-1}|\leq |a|$.
Hence, by Proposition \ref{401b}, $p=r\cdot q^t$
for some $\alpha$-peak $r$ with $|r|\leq 2|pqp^{-1}|^2\leq 2|a|^2$.
Similarly, $p=r'\cdot (q')^{t'}$
for some $\alpha$-peak $r'$ with $|r'|\leq 2|a'|^2\leq 2|a|^2$.
Now
\de{15ma15z}{
\ac{p}\geq t=(|p|-|r|)/|q|\geq (8|a|^3-2|a|^2)/|a|\geq 6|a|^2\geq 3\max\{|r|,|r'|,|q|,|q'|\},
}
Hence by \rf{28ja15ax}, $q=d^n$ for some $d$ and $n\geq 0$, with
$|d|\leq|q'|\leq|a'|$.
As $t\geq 1$ by \rf{15ma15z}, $q$ is cyclically reduced, hence also $d$ is cyclically reduced.
As $q$ is connected, also $d$ is connected.

Let $c$ be a conjugate of $a$.
Then $c$ is a conjugate of $p^{-1}ap=qu$, with $u$ independent of $q$.
Hence by \rf{15ma15x}, $c$ contains a segment $c'$ which is a conjugate of $q=d^n$.
Now $n=|q|/|d|\geq \ac{p^{-1}ap}/|a'|>2|a'|\geq 2|d|$.
Hence, by Proposition \ref{1fe15c}, $c'$ contains $s:=d^{n-2|d|}$ as segment.
Then $qs^{-1}=d^{2|d|}$, hence $\ac{p^{-1}aps^{-1}}=\ac{qs^{-1}}\leq |qs^{-1}|
=2|d|^2\leq 2|a'|^2$.
\bx 

\subsectz{The function $\mu_{a,H}:x\mapsto\min(a^{-1}x^{\uparrow}\hspace*{-4.5pt}H)$}

The following function $\mu_{a,H}:G\to G$ forms an important ingredient in our algorithm.
Fixing $a\in G$ and a closed set $H\subseteq G$, it is defined by
\dyyz{
\mu_{a,H}(x):=\min(a^{-1}x^{\uparrow}H)
}
for $x\in G$, which is well-defined as $a^{-1}xH$ is convex by \rf{p7f}.
So for each $x\in G$:
\de{26no14a}{
\mu_{a,H}(x)\leq a^{-1}x.
}
Moreover, for all $x,y\in G$:
\dy{9ma15a}{
$\mu_{a,H}(x)\leq y$ if and only if $a\in x^{\uparrow}H(y^{\downarrow})^{-1}$,
}
since $\min(a^{-1}x^{\uparrow}H)\leq y$ if and only if
$a^{-1}bh=c$ for some $b\in x^{\uparrow}$, $h\in H$ and $c\in y^{\downarrow}$.

\prop{11no14a}{
Let $a\in G$ and $H\subseteq G$ be closed, and set $\mu:=\mu_{a,H}$.
Then for all $x,y\in G$:%
\di{25ok14i}{
\nr{i} if $x\leq y$ then $\mu(x)\leq\mu(y)$,
\nrs{ii}
$\mu(x\wedge y)\leq\mu(x)\wedge \mu(y)$,
\nrs{iii}
if $x\vee y$ is finite, then
$\mu(x\vee y)=\mu(x)\vee\mu(y)$.
}
}

\pf
Since $x\leq y$ implies $x^{\uparrow}\supseteq y^{\uparrow}$, we have (i).
Then (ii) follows from (i).
To see (iii), we have $\mu(x)\vee\mu(y)\leq\mu(x\vee y)$ by (i).
In particular, $\mu(x)\vee\mu(y)$ is finite.
To see the reverse inequality,
set $d:=\mu(x)=\min(a^{-1}x^{\uparrow}H)$ and $e:=\mu(y)=\min(a^{-1}y^{\uparrow}H)$.
Then, by \rf{9ma15a} and as both $d^{-1}$ and $e^{-1}$ belong to $((d\vee e)^{\downarrow})^{-1}$,
\dyyz{
a\in
x^{\uparrow}Hd^{-1}
\cap
y^{\uparrow}He^{-1}
\subseteq
x^{\uparrow}H((d\vee e)^{\downarrow})^{-1}
\cap
y^{\uparrow}H((d\vee e)^{\downarrow})^{-1}
=
(x\vee y)^{\uparrow}H((d\vee e)^{\downarrow})^{-1},
}
where the equality follows from \rf{p1a}
(as $x^{\uparrow}\cap y^{\uparrow}=(x\vee y)^{\uparrow}$).
So by \rf{9ma15a}, $\mu(x\vee y)=\min(a^{-1}(x\vee y)^{\uparrow}H)\leq d\vee e$.
\bx

The composition of functions $\mu_{a,H}$ have the following property.
Let $x,a,a'\in G$ and let $H$ and $H'$ be closed.
Then
\de{12no14a}{
\mu_{a',H'}(\mu_{a,H}(x))\leq \mu_{aa',HH'}(x).
}
Indeed, by the definitions of $HH'$ and $\mu_{aa',HH'}(x)$,
there exist $x'\geq x$ and $c\in H$, $c'\in H'$ with $\mu_{aa',HH'}(x)=(aa')^{-1}x'cc'$.
Now $\mu_{a,H}(x)\leq a^{-1}x'c$.
Hence, using \rf{25ok14i}(i),
\dyyz{
\mu_{a',H'}\circ\mu_{a,H}(x)
\leq
\mu_{a',H'}(a^{-1}x'c)
\leq
(a')^{-1}(a^{-1}x'c)c'
=
\mu_{aa',HH'}(x).
}

\subsectz{Polynomial-time algorithms}

Let $I$ be an ideal and $x\in G$.
Then there is a unique largest element $y\leq x$ with $y\in I$.
We can find it in polynomial time if membership of $I$ can be tested in polynomial time:
\dyz{
Let $I$ be an ideal of which we can test membership in polynomial time.
Then for any $x\in G$, we can find the maximal element $y\leq x$ with $y\in I$ in polynomial time.
}
To see this, grow a word $y\leq x$ with $y\in I$, starting with $y=1$.
If there is a minimal symbol in $y^{-1}x$ with $y\alpha\in I$, replace $y$ by $y\alpha$.
If no such $\alpha$ exists, $y$ is as required.

Note that $y$ is the closest (with respect to $\dist$) element in $I$ to $x$.
Hence, by the left-invariance of the distance function,
$x^{-1}y$ is the closest element in $x^{-1}I$ to 1.
That is: $x^{-1}y=\min(x^{-1}I)$.
Therefore:
\dy{1no14f}{
Let $I$ be an ideal of which we can test membership in polynomial time.
Then for any $z\in G$, we can find $\min(zI)$ in polynomial time.
}

Note that we can test membership of $x^{\uparrow}$ and of $x^{\downarrow}$ in polynomial time.
Hence:
\dy{1no14g}{
For any $y,x\in G$, we can find $\min(y^{-1}x^{\uparrow})$ in polynomial time.
}
This follows from \rf{1no14f} setting $I:=x^{-1}x^{\uparrow}$ and $z:=y^{-1}x$.

If $H$ is closed, then for any $y\in G$, $y^{-1}x^{\uparrow}H$ is convex
(by \rf{p7f}).
Its minimum  can be found in polynomial time:
\dy{23ok14b}{
Let $H$ be a closed set of which we can test membership in polynomial time.
Then for any $x,a\in G$,
$\mu_{a,H}(x)=\min(a^{-1}x^{\uparrow}H)$ can be found in polynomial time.
}
Indeed, by \rf{p7f}, $\min(a^{-1}x^{\uparrow}H)=\min(\min(a^{-1}x^{\uparrow})H)$.
Hence \rf{23ok14b} follows from \rf{1no14g} and \rf{1no14f}.

Finally,
\dy{1no14h}{
Let $H$ be a closed set of which we can test membership in polynomial time.
Then for any $x,y,z\in G$, we can test in polynomial time if
$y$ belongs to $x^{\uparrow}H(z^{\uparrow})^{-1}$.
}
Indeed, $y\in x^{\uparrow}H(z^{\uparrow})^{-1}$ if and only if
$y^{-1}x^{\uparrow}H\cap z^{\uparrow}\neq\emptyset$.
This is the case if and only if $z^{-1}y^{-1}x^{\uparrow}H\cap z^{-1}z^{\uparrow}\neq\emptyset$.
The latter statement is equivalent to:
$\min(z^{-1}y^{-1}x^{\uparrow}H)$ belongs to the ideal $z^{-1}z^{\uparrow}$.
As this minimum can be determined in polynomial time by \rf{23ok14b} and as membership of $z^{-1}z^{\uparrow}$
can be tested in polynomial time, we have proved \rf{1no14h}.

\sect{13me03b}{The cohomology feasibility problem}

Let $D=(V,E)$ be a directed graph and let $G$ be a group.
Two functions $\phi,\psi:E\to G$ are called {\em cohomologous} if there exists a
function $f:V\to G$ such that $\psi(e)=f(u)^{-1}\phi(e)f(w)$ for each edge $e=(u,w)$.
One directly checks that this gives an equivalence relation.

We give a polynomial-time algorithm for the following {\em cohomology feasibility problem}
for graph groups:
\di{a3}{
\item[\em given:] a directed graph $D=(V,E)$, an undirected graph $([k],F)$, a function
$\phi:E\to G_F$, and for each edge $e$ a closed set $H(e)\subseteq G_F$,
\items{\em find:} a function $\psi:E\to G_F$ such that $\psi$ is cohomologous to $\phi$
and such that $\psi(e)\in H(e)$ for each $e\in E$.
}

The running time of the algorithm for this problem is bounded by a polynomial
in $n:=|V|$, $m:=|E|$, $\sigma:=\max\{|\phi(e)|\mid e\in E\}$, and the maximum time needed to test
if any word of polynomial size belongs to $H(e)$ (over all edges $e$).
(The number $k$ of generators can be bounded by $m\sigma$, since we may assume that all
generators occur among the $\phi(e)$.)
More precisely, there exist polynomials $p_1$ and $p_2$ such that the problem takes time
$p_1(n,m,\sigma)\tau_H(p_2(n,m,\sigma))$,
where $\tau_H(x)$ is the time needed to test membership of the $H(e)$ for words
of size at most $x$.

Note that, by the definition of cohomologous, equivalent to finding a function $\psi$ as in \rf{a3},
is finding a function $f:V\to G_F$ satisfying:
\dy{a3a}{
$f(u)^{-1}\phi(e)f(w)\in H(e)$ for each edge $e=(u,w)$.
}
We call such a function $f$ {\em feasible}.

We can assume that
\dy{2no14a}{
$|\phi(e)|\leq 1$ for each edge $e$.
}
Indeed, if $\phi(e)=xy$ for edge $e=(u,w)$, we can split the edge into two edges $(u,v),(v,w)$, where
$v$ is a new vertex, and define $\phi(u,v):=x$, $\phi(v,w):=y$, $H(u,v):=H(e)$, and $H(v,w):=\{1\}$.
The new problem is equivalent to the original problem: if $f$ is a solution to the original problem,
we can set $f(v):=yf(w)$, and obtain a solution for the new problem; conversely, if $f$ is a solution to
the original problem, forgetting the value of $f$ on $v$, we obtain a solution to the original problem.

\subsect{s3a}{Pre-feasible functions}

Given the input of the cohomology feasibility problem \rf{a3}, we
call a function $f:V\to G_F$ {\em pre-feasible} if for each edge
$e=(u,w)$ of $D$ there exist $x\geq f(u)$ and $z\leq f(w)$ such that $x^{-1}\phi(e)z\in H(e)$.
Clearly, each feasible function is pre-feasible.
There is a trivial pre-feasible function $f$, defined by $f(v):=1$ for each $v\in V$.
Note that $f$ is pre-feasible if and only if
\dy{25ok14j}{
$\mu_{\phi(e),H(e)}(f(u))\leq f(w)$ for each edge $e=(u,w)$.
}

The collection $G_F^V$ of all functions $f:V\to G_F$ can be partially ordered by: $f\leq g$ if
and only if $f(v)\leq g(v)$ for each $v\in V$.
Then $G_F^V$ forms a lattice if we add an element $\infty$ at infinity.
Let $\wedge$ and $\vee$ denote meet and join.
Then \rf{25ok14i}(ii) and (iii) directly give:
\dy{P25a}{\em
Let $f_1$ and $f_2$ be pre-feasible functions.
Then $f_1\wedge f_2$ and, if $f_1\vee f_2<\infty$, $f_1\vee f_2$ are pre-feasible again.
}

It follows that for each function $f:V\to G_F$ there is a unique smallest pre-feasible function
$\bar{f}\geq f$, provided that there exists at least one pre-feasible function $g\geq f$.
If no such $g$ exists we set $\bar{f}:=\infty$.
By \rf{P25a}, $\overline{f\vee g}=\bar{f}\vee\bar{g}$ for any two functions $f,g$ with
$f\vee g$ finite.

\subsect{s4}{A subroutine finding $\bar{f}$}

Condition \rf{25ok14j} suggests a `constraint satisfaction' algorithm to find $\bar{f}$ for a given
function $f$.
Let input $D=(V,E),F,\phi,H$ for the cohomology feasibility problem be given.
For any edge $e$, we can determine $\mu_{\phi(e),H(e)}(f(u))$ in polynomial time, by \rf{23ok14b}.

\vspace{4mm}
\noindent
{\bf Subroutine to find $\bar{f}$:} 
Find an edge $e=(u,w)$ for which $m_e:=\mu_{\phi(e),H(e)}(f(u))\not\leq f(w)$.
If $m_e\vee f(w)$ is finite, reset $f(w):=m_e\vee f(w)$ and iterate.
If $m_e\vee f(w)=\infty$, output $\bar{f}:=\infty$.
If no such edge $e$ exists, output $\bar{f}:=f$.

\medskip
Then the output of the subroutine (if any) is correct.
For let $f'$ be the reset function.
If $\bar{f}$ is finite, then $f\leq f'\leq\bar{f}$, since
$f'(w)=\mu_{\phi(e),H(e)}(f(u))\vee f(w)
\leq
\mu_{\phi(e),H(e)}(f(u))\vee \bar{f}(w)
=\bar{f}(w)$,
since $\bar{f}$ is pre-feasible.
So in this case $\bar{f'}=\bar{f}$.
This moreover implies that if $m_e\vee f(w)=\infty$ then $\bar{f}=\infty$.

\subsect{s5a}{Running time of the subroutine}

For each walk $P=e_1e_2\ldots e_t$, where $e_1,\ldots,e_t$ are consecutive edges of $D$,
we set $\phi(P):=\phi(e_1)\phi(e_2)\ldots\phi(e_t)$ and $H(P):=H(e_1)H(e_2)\linebreak[1]\ldots\linebreak[1] H(e_t)$.

We will study the running time of the subroutine under the condition that the cohomology feasibility
problem has a solution, or more weakly, that
\dy{15de14b}{
for each closed walk $C$, $H(C)$ contains a conjugate of $\phi(C)$.
}

For any function $f:V\to G_F$, let $|f|:=\max\{|f(v)|\mid v\in V\}$.

\prop{21ok14e}{
If \rf{2no14a} and \rf{15de14b} hold, then for any $f:V\to G_F$ with $\bar f$ finite:
\dez{
|\bar f|\leq k|f|+220k(2nk)^9.
}
}

\pf
Define $c:=2nk$ and $m:=3c^2$, and suppose $|\bar{f}|>k|f|+220kc^9$.
Consider the (infinite) directed graph $\DD$ with vertex set $\VV:=\{(v,p)\mid v\in V$, $p$ peak$\}$ and edge set
all pairs $((v,p),(w,q))\in\VV\times\VV$ with $e=(v,w)\in E$ and $q\leq\mu_{\phi(e),H(e)}(p)$.
Let $\SS$ be the set of all vertices $(v,p)$ of $\DD$ with $p\leq f(v)$.
Then for each $v\in V$ and peak $p$:
\dyz{
$\bar{f}(v)=\bigvee\{p\mid$ there exists a walk in $\DD$ from $\SS$ to $(v,p)\}$.
}
This follows with \rf{16no14a} and Proposition \ref{11no14a}(iii).
Hence:
\dy{23no14b}{
if $p\leq \bar{f}(v)$, then there exists a walk in $\DD$ from $\SS$ to $(v,p)$.
}

Since $|\bar{f}|>k|f|+220kc^9$, there exists $w\in V$ with $|\bar{f}(w)|>k|f|+220kc^9$.
Hence there exists a peak $q\leq\bar{f}(w)$ with $|q|>|f|+220c^9$
(since $\bar{f}(w)$ is a join of at most $k$ peaks).
By \rf{23no14b}, there is a walk $\PP$ in $\DD$ from $\SS$ to $(w,q)$.
Choose a shortest such walk $\PP$; let it have length $\ell$.
Since $|\phi(e)|\leq 1$ by assumption,
$|p|\leq|\phi(e)^{-1}p'|\leq|p'|+1$ for each edge $((v,p),(v',p'))$ of $\DD$, where $e=(v,v')$.
Hence $\ell\geq|q|-|f|>220c^9$.

Let $\PP$ traverse vertices $(v_0,q_0),(v_1,q_1)\ldots,(v_{\ell},q_{\ell})$ of $\DD$, in this order.
So $v_{\ell}=w$ and $q_{\ell}=q$.
For each $u\in V$ and each symbol $\alpha$ let $I_{u,\alpha}$ denote the set of indices
$j\in\{\ell-3c^3,\ldots,\ell\}$
such that $v_j=u$ and $q_j$ is an $\alpha$-peak.
Then there exist $u\in V$ and a symbol $\alpha$ such that $|I_{u,\alpha}|>3c^3/c=m$.
Choose $j_0<j_1<\cdots<j_m$ in $I_{u,\alpha}$.
Set $p_i:=q_{j_i}$ for $i=0,1,\ldots,m$.
So each $p_i$ is an $\alpha$-peak and $p_0<p_1<\cdots<p_m$.
We will apply Proposition \ref{26no14b} to $p:=p_0$.

Let $C_i$ be the $u-u$ walk $v_{j_{i-1}},v_{j_{i-1}+1},\ldots,v_{j_i-1},v_{j_i}$ in $D$, and
let $C:=C_1C_2\ldots C_m$ (the concatenation of $C_1,\ldots,C_m$).
Then $|C|\leq 3c^3$ and so there exists $i\in\{1,\ldots,m\}$ with $|C_i|\leq |C|/m\leq 3c^3/m=c$.

Let $a:=\phi(C)^{-1}$ and $a'':=\phi(C_i)^{-1}$.
Then $|a|=|\phi(C)|\leq|C|\leq 3c^3$ and $|a''|=|\phi(C_i)|\leq|C_i|\leq c$.
This gives
\dez{
|p_0|\geq |q|-3c^3\geq 220c^9-3c^3\geq 216c^9=8(3c^3)^3\geq 8|a|^3.
}
Moreover, $p_0\leq p_m\leq\mu_{\phi(C),H(C)}(p_0)\leq ap_0$ and similarly $p_i\leq p_{i+1}\leq a''p_i$.
As $\alpha$ occurs in $p_i^{-1}p_{i+1}$, $\alpha$ also occurs in $p_i^{-1}a''p_i$.
By \rf{15ma15f}, $p_0\leq a'p_0$ for some $a'$ with $|a'|\leq |a''|$ and
$\ac{p_0^{-1}a'p_0}=\ac{p_i^{-1}a''p_i}$.

As we assume that \rf{15de14b} holds, some conjugate $x$ of $\phi(C)$ belongs to $H(C)$.
So $x^{-1}$ is a conjugate of $a$.
Hence, by Proposition \ref{26no14b},
$x^{-1}$ has a segment $s$ such that $\ac{p_0^{-1}ap_0s^{-1}}\leq 2|a'|^2\leq 2c^2<3c^2=m$.
Now $p_m\leq \min(\phi(C)^{-1}p_0^{\uparrow}H(C))\leq ap_0s^{-1}$, since $s^{-1}\in H(C)$, as $H(C)$ is closed and $x\in H(C)$.
So $\ac{p_0^{-1}p_m}< m$, contradicting the fact that $p_m$ contains at least
$m$ $\alpha$'s more than $p_0$.
\bx 

This implies, where $n:=|V(D)|$,
and where, for any $\sigma\in\oZ_+$, $\tau_{H}(\sigma)$ is the maximum time needed to test if any word of
size $\leq\sigma$ belongs to $H(e)$, for any given edge $e$.

\thm{t22a}{
There exist polynomials $p_1$ and $p_2$ such that, if \rf{2no14a} and \rf{15de14b} hold, then
the running time of the subroutine is bounded by $p_1(n,k,|f|)\tau_H(p_2(n,k,|f|))$.
}

\pf 
At each iteration, we increase $|f(v)|$ for some vertex $v$.
Hence Proposition \ref{21ok14e} implies that, if $\bar{f}$ is finite,
the number of iterations is bounded by some polynomial $p_1$ in $n$, $k$, $|f|$.
If the subroutine exceeds this number of iterations, we conclude that$\bar{f}=\infty$.

Since in each iteration, the reset $f'$ satisfies $|f'|\leq|\bar f|$, and since $|\bar f|$ is bounded by
a polynomial $p_2$ in $n$, $k$, $|f|$, in each iteration we only need to test membership of words
of size at most $p_2(n,k,|f|)$.
\bx

\subsect{s5}{A polynomial-time algorithm for the cohomology feasibility problem for graph groups}

We now describe the algorithm for the cohomology feasibility problem for graph groups.
Let input $D=(V,E)$, $F$, $\phi$, $H$ of \rf{a3} be given.

Let $\FF$ be the collection of all functions $f:V\to G_F$ such that for each edge $e=(u,w)$ of $D$ there
exist $x\geq f(u)$ and $z\geq f(w)$ satisfying $x^{-1}\phi(e)z\in H(e)$;
equivalently:
\dez{
\phi(e)\in f(u)^{\uparrow}H(e)(f(w)^{\uparrow})^{-1}.
}
This can be tested in polynomial time by \rf{1no14h}.
So for any given function $f$ one can check in polynomial time whether $f$ belongs to $\FF$.
Trivially, if $f\in\FF$ and $g\leq f$ then $g\in\FF$.
Moreover:

\dy{p14A}{
Let $f_1,\ldots,f_t$ be functions such that $f_i\vee f_j\in\FF$ for all $i,j$.
Then $f:=f_1\vee\cdots\vee f_t\in\FF$.
}

\pf
We must show that for each edge $e=(u,w)$, $\phi(e)$ belongs to
${f(u)}^{\uparrow}H(e)({f(w)}^{\uparrow})^{-1}$.
Since $f_i\vee f_j\in\FF$ for all $i,j$, we know
\de{d14D}{
\phi(e)\in f_i(u)^{\uparrow}H(e)f_j(w)^{\uparrow-1}}
for all $i,j$.
Hence by \rf{p1a},
\dyy{d14E}{
\phi(e)\in\bigcap_i\bigcap_jf_i(u)^{\uparrow}H(e)(f_j(w)^{\uparrow})^{-1}
=\big(\bigcap_if_i(u)^{\uparrow}\big)H(e)\big(\bigcap_j(f_j(w)^{\uparrow})^{-1}\big)=
{f(u)}^{\uparrow}H(e)({f(w)}^{\uparrow})^{-1}.
}
Here $i$ and $j$ range over $1,\ldots,t$.
\bx

In the following theorem, `solvable in polynomial time' means as before that there exist
polynomials $p_1$ and $p_2$ such that the problem is solvable in time
$p_1(n+m,k,\rho)\tau_H(p_2(n+m,k,\rho))$, where
$n:=|V(D)|$, $m:=|E(D)|$, $\rho$ is the maximum of $|\phi(e)|$ over all $e\in E$, and where
$\tau_H(\sigma)$ again is the maximum time needed to test if any word of
size $\leq\sigma$ belongs to $H(e)$, for any given edge $e$.

\thm{t29a}{
The cohomology feasibility problem for graph groups is solvable in polynomial time.
}

\pf
We can assume again that $|\phi(e)|\leq 1$ for each edge $e$.
Moreover, we can assume that with each edge $e=(u,w)$ also $e^{-1}=(w,u)$ is an edge,
with $\phi(e^{-1})=\phi(e)^{-1}$ and $H(e^{-1})=H(e)^{-1}$.

For any $e=(u,w)\in E$, let $f_e$ be the function defined by
\de{d14B}{
f_e(v):=
\begin{cases}
\phi(e) & \text{ if $v=u$,} \\
1& \text{ if $v\neq u$.}
\end{cases}
}
Let $L$ be the set of pairs $\{e,e^{-1}\}$ from $E$ such that $\phi(e)\not\in H(e)$.
Let $N$ be the collection of all pairs $\{e,d\}$ from $E$ such that the function
$\bar{f}_e\vee\bar{f}_d=\infty$, or is finite and does not belong to $\FF$ (possibly $e=d$).

Choose a subset $B$ of $E$ such that $B$ intersects each pair in $L$ and such that $B$ contains no
pair in $N$.
This is a special case of the 2-satisfiability problem, and hence can be solved in polynomial time.
Assuming that there exists a feasible function $f$, then $B$ exists,
as $B:=\{e=(u,v)\in E\mid \phi(e)\leq f(u)\}$ would have the required properties.

If we find $B$, define $f$ by:
\de{d12N}{
f(v):=\bigvee_{e\in B}\bar{f}_e.
}
We are done by proving that $f$ is feasible.
Since $\bar{f}_e\vee\bar{f}_b<\infty$ for each pair $\{e,d\}\subseteq B$, we know $f<\infty$.
Moreover, $f$ is the join of a finite number of pre-feasible functions, and hence $f$ is pre-feasible.
So by \rf{15de14a} it suffices to prove that for each edge $e=(u,w)$: 
\di{d14C}{
\nr{i} there exist $x\geq f(u)$ and $z\geq f(w)$ such that $x^{-1}\phi(e)z\in H(e)$,
\nrs{ii} there exist $x\leq f(u)$ and $z\leq f(w)$ such that $x^{-1}\phi(e)z\in H(e)$.
}

To prove \rf{d14C}(i), note that it is equivalent to: $f\in\FF$.
As $\bar{f}_e\vee\bar{f}_b\in\FF$ for all
$a,b\in B$, \rf{p14A} gives $f\in\FF$.

To prove \rf{d14C}(ii), if it does not hold then
$\phi(e)\not\in H(e)$, hence $\{e,e^{-1}\}\in L$.
So $e$ or $e^{-1}$ belongs to $B$.
By symmetry, we can assume that $e\in B$.
So $f_e\leq f$, and therefore $\phi(e)\leq f(u)$.
So we can take $x:=\phi(e)$ and $z:=1$ in \rf{d14C}(ii).
\bx

An analysis of this algorithm shows that the cohomology feasibility problem has a solution
if and only if for each vertex $u$ and each pair $C,C'$ of (undirected) $u-u$ walks in $D$ there exists
$x\in G$ such that $x^{-1}\phi(C)x\in H(C)$ and $x^{-1}\phi(C')x\in H(C')$.
This condition is trivially necessary.

\sect{2de14c}{Planar graphs}

We repeat the {\em partially disjoint paths problem} for directed planar graphs:
\di{27no14a}{
\item[{\em given:}]
a directed planar graph $D=(V,E)$, vertices $r_1,s_1,\ldots,r_k,s_k$ of $D$, and a set $F$ of pairs $\{i,j\}$ with $i,j\in [k]$,
\items{\em find:}
a $k$-tuple $\PP=(P_1,\ldots,P_k)$, where $P_i$ is a directed $r_i-s_i$ path $P_i$,
for $i=1,\ldots,k$, such that $P_i$ and $P_j$ are disjoint whenever $\{i,j\}\in F$.
}

We can assume without loss of generality:
\dy{26de14a}{
$r_1,s_1,\ldots,r_k,s_k$ are distinct, each
$r_i$ has outdegree 1 and indegree 0, and each $s_i$ has indegree 1 and outdegree 0.
}

Again, let $G_F$ be the graph group generated by $g_1,\ldots,g_k$ and relations $g_ig_j=g_jg_i$ whenever
$\{i,j\}\not\in F$.
For each solution $\PP$ of \rf{27no14a}, let $\chi_{\PP}:E\to G_F$ be defined by:
\dez{
\chi_{\PP}(e):=\prod_{i\atop P_i\text{ traverses }e}g_i.
}
The order in which we take this product is irrelevant, since if both $P_i$ and $P_j$ traverse $e$,
then $\{i,j\}\not\in F$ and hence $g_ig_j=g_jg_i$.

Let $\FF$ be the collection of faces of $D$.
Call $\phi,\psi:E\to G_F$ {\em homologous} if there exists $f:\FF\to G_F$
such that for each edge $e$: $\psi(e)=f(F)^{-1}\phi(e)f(F')$, where $F$ and $F'$
are the left-hand and the right-hand face at $e$, respectively
(seen when traversing $e$ in forward direction).

\subsectz{Finding partially disjoint paths of prescribed homology}

We first consider the homology version of the partially disjoint paths problem:
\di{27no14c}{
\item[{\em given:}]
a directed planar graph $D=(V,E)$,
vertices $r_1,s_1,\ldots,r_k,s_k$ of $D$ satisfying \rf{26de14a},
a set $F$ of pairs $\{i,j\}$ with $i,j\in [k]$,
and a function $\phi:E\to G_F$,
\items{\em find:}
a solution $\PP$ of \rf{27no14a} such that $\chi_{\PP}$ is homologous to $\phi$.
}

\prop{28no14b}{
Problem \rf{27no14c} is solvable in polynomial time.
}

\pf
We can assume that problem \rf{27no14c} has a solution --- that is, $\phi$ is homologous to
$\chi_{\PP}$ for some solution $\PP$ of \rf{27no14a}.

Let $\FF$ be the collection of faces of $D$.
Consider the dual directed graph $D^*=(\FF,E^*)$, where for each edge $e$ of $D$ there is a directed
edge $e^*\in E^*$ from the face at the left-hand side of $e$ to the face at the right-hand side of $e$.
We define $\hat\phi(e^*):=\phi(e)$ and
\dy{28no14d}{
$H(e^*):=\{\prod_{i\in I}g_i\mid I\subseteq[k], I$ stable set in $([k],F)\}$,
}
where $I$ is {\em stable} if it contains no pair in $F$ as subset.
Note that $H(e^*)$ is a closed subset of $G_F$.

We extend the planar graph $D^*$ by a number of further `nonplanar' edges, as follows.
Consider any vertex $v\not\in\{r_1,s_1,\ldots,r_k,s_k\}$ of $D$ and two faces $F$ and $F'$ of $D$ incident with
$v$.
Let $e_1,\ldots,e_t$ be the edges incident with $v$ that are crossed when going clockwise from
$F$ to $F'$ around $v$.
Then add to $D^*$ an edge $e_{v,F,F'}$ from $F$ to $F'$, and define
\dez{
\hat\phi(e_{v,F,F'}):=\phi(e_1)^{\sigma_1}\ldots\phi(e_t)^{\sigma_t},
}
where, for $j=1,\ldots,t$, $\sigma_j:=1$ if $e_j$ is oriented
away from $v$ and $\sigma_j:=-1$ if $e_j$ is oriented towards $v$.
Note that, as by assumption $\phi$ is homologous to $\chi_{\PP}$ for some solution $\PP$ of \rf{27no14a},
we necessarily have $\hat\phi(e_{v,F',F})=\hat\phi(e_{v,F,F'})^{-1}$.

Moreover, define
\dyz{
$H(e_{v,F,F'}):=\{\prod_{i\in I}g_i^{\tau(i)}\mid I\subseteq[k], I$ stable set in $([k],F),
\tau:I\to\{+1,-1\}\}$.
}
Also $H(e_{v,F,F'})$ is a closed subset of $G_F$.

Let $\hat D$ be the extended directed graph, and
consider the cohomology feasibility problem $\Pi$ for $(\hat D,\hat\phi)$, in which we require that
the output is only weakly allowed on the edges in $E'$.
As, by assumption, $\phi$ is homologous to $\chi_{\PP}$ for some $\PP$,
problem $\Pi$ has a solution, namely $\chi_{\PP}$.
Conversely, let $\psi:E(\hat D)\to G_F$ be any solution of $\Pi$.
Define $\check\psi:E(D)\to G_F$ by $\check\psi(e):=\psi(e^*)$ for $e\in E(D)$.
Then $\check\psi$ is equal to $\chi_{\PP}$ for some $\PP$.
This because, by our assumption, $\phi$ is homologous to $\chi_{\PP}$ for some $\PP$.
Hence, for the edge $e$ incident with $r_i$, $\phi(e)$ is conjugate to $g_i$ (as by \rf{26de14a} $e$ is incident
at both sides to the same face), therefore
$\check\psi(e)$ is conjugate to $g_i$.
Since $\check\psi(e)$ is allowed, it follows that in fact $\check\psi(e)=g_i$.
So only $P_i$ traverses $e$, and in the forward direction.
Similarly for the edge incident with $s_i$.

Therefore, the proposition follows from Theorem \ref{t29a}.
\bx

\subsect{28no14e}{Enumerating homologies of disjoint paths}

We finally describe an algorithm that finds, for any input of \rf{27no14a},
a collection $\Phi$ of functions $\phi:E\to G_F$ with the property that
\dy{29no14a}{
for each solution $\PP$ of \rf{27no14a}, $\chi_{\PP}$ is homologous to at least one $\phi\in\Phi$.
}
So, although there exist infinitely many homology classes (if $F\neq\emptyset$),
in our algorithm we can restrict ourselves to a number of homology classes that, for fixed $k$,
is bounded by a polynomial in the size of the graph.

\prop{28no14c}{
Fixing $k$, a collection $\Phi$ satisfying \rf{29no14a} can be found in polynomial time.
}

\pf
Again, we can assume \rf{26de14a}.
Moreover, we can assume that $D$ is weakly connected and that (for the convenience of the exposition)
each $i\in[k]$ is contained in at least one
pair $\{i,j\}$ in $F$ (otherwise we can easily reduce the problem).
We also can assume that each vertex $v\neq r_1,s_1,\ldots,r_k,s_k$
has total degree $\deg(v)$ equal to 3: 
replace $v$ by a directed circuit of length $\deg(v)$ and attach the edges incident with $v$
to the different vertices of the circuit (in a planar manner of course).
Any $\phi$ found for the modified graph can be `shrunk' to the smaller graph.

Choose a spanning tree $T$ in $D$.
We will consider graphs $T'$ obtained from $T$ by replacing each edge $e$ by a number (possibly 0) of parallel edges.
These edges form a {\em parallel class}, denoted by $\pi_e$.
Each such graph $T'$ is trivially planar, by drawing the edges properly parallel in the plane.

We moreover consider undirected walks in such graphs $T'$.
(An {\em undirected walk} may traverse edges in any direction.)
Call undirected walks $W$ and $W'$ in $T'$ {\em crossing} if there is a vertex $v$ and distinct edges
$e_1,e_2,e_3,e_4$ of $T'$ incident with $v$, in clockwise or counterclockwise order,
such that $W$ traverses $e_1$ and $e_3$ consecutively, and $W'$ traverses $e_2$ and $e_4$ consecutively.
If $W=W'$, we say that $W$ is {\em self-crossing}.

In particular, we consider $k$-tuples $\WW=(W_1,\ldots,W_k)$ of undirected walks in $T'$ such that
\di{27no14b}{
\nr{i} $W_i$ runs from $r_i$ to $s_i$ and is not self-crossing, for each $i=1,\ldots,k$,
\nrs{ii} $W_i$ and $W_j$ are not crossing, for each $\{i,j\}\in F$,
\nrs{iii} each edge of $T'$ is traversed by precisely one $W_i$.
}
The last condition implies that $T'$ is determined by $W_1,\ldots,W_k$.

For each $k$-tuple $\WW$ satisfying \rf{27no14b}, define $\phi_{\WW}:E\to G_F$ as follows.
If $e$ is an edge of $D$ not in $T$, set $\phi_{\WW}(e):=1$.
If $e=(u,w)$ is an edge of $T$, let
$e_1,\ldots,e_t$ be the edges in $\pi_e$, from left to right with respect to the orientation
$(u,w)$ of $e$.
Let $\alpha_j:=g_i$ if $e_j$ is traversed by $W_i$ in the direction from $u$ to $w$, and
let $\alpha_j:=g_i^{-1}$ if $e_j$ is traversed by $W_i$ in the direction from $w$ to $u$.
Define $\phi_{\WW}(e):=\alpha_1\ldots\alpha_t$.
Then:
\dy{2de14a}{
For each solution $\PP$ of \rf{27no14a}, there exists $T'$ and
$\WW=(W_1,\ldots,W_k)$ satisfying \rf{27no14b}
such that $\chi_{\PP}$ and $\phi_{\WW}$ are homologous and such that,
for each $i=1,\ldots,k$, each parallel class in $T'$ is traversed at most $2|E|$ times by $W_i$.
}
To see this, reroute $P_1,\ldots,P_k$ along $T$ as follows.
For each $e$ in $E(D)\setminus E(T)$, let $Q_e$ be the path in $T$ connecting the ends of $e$.
Order the edges in $E(D)\setminus E(T)$ as $e_1,e_2,\ldots,e_m$ such that if $Q_{e_j}$ is longer than $Q_{e_i}$ then $j>i$.
Then for $j=1,\ldots,m$, if $P_i$ traverses $e_j$, reroute $P_i$ along $Q_{e_j}$; that is,
add edges parallel to the edges in $Q_{e_j}$, in the disk enclosed by $e_j$ and $Q_{e_j}$,
and replace $e_j$ in $P_i$ by the new edges (in order).
This gives $T'$ and $W_1,\ldots,W_k$ as required, proving \rf{2de14a}.

So to cover all homology classes of solutions of problem \rf{27no14a}, it suffices to enumerate all $T'$
and $W_1,\ldots,W_k$ satisfying \rf{27no14b}.

In fact, we can assume that each $W_i$ is {\em non-returning} in the following sense.
Let $W_i$ traverse edges $e$, vertex $v$, and edges $e'$ consecutively.
\dy{2fe15a}{
(i) If $v\not\in\{r_1,s_1,\ldots,r_k,s_k\}$, then $e$ and $e'$ belong to different parallel classes
incident with $v$.\\
(ii) If $v\in \{r_j,s_j\}$ for some $j\in[k]$, then $e$ and $e'$ enclose the starting or ending edge of $W_j$.
}
This can be attained as follows.
Suppose $W_i$, $e$, $v$, $e'$ violate \rf{2fe15a}.
Fixing $v$, choose $W_i,e,e'$ such that the number of edges inbetween of (that is, enclosed by)
$e$ and $e'$ is as small as possible.
Then each edge inbetween of $e$ and $e'$ is traversed by some $W_j$ with $j\neq i$ (as $W_i$ is not
self-crossing) and $\{i,j\}\not\in F$ (as $W_i$ and $W_j$ cross).
So deleting $e$ and $e'$ from $W_i$ and from $T'$, gives a walk system $\WW'$ again
satisfying \rf{27no14b}, with $\phi_{\WW'}=\phi_{\WW}$, and with a smaller total length.
Iterating this, we end up with each $W_i$ non-returning.

This implies that if $v\in V\setminus\{r_1,s_1,\ldots,r_k,s_k\}$ has degree 1 in $T$, it is incident with
no edges in $T'$.
Delete such vertices from $T$ repeatedly.
Let $T_0$ be the final graph.
It is a tree with maximum degree 3 and with $2k$ vertices of degree 1 (namely, $r_1,s_1,\ldots,r_k,s_k$).
Hence $T_0$ has $2k-2$ vertices of degree 3.
The vertices of degree 1 and 3 are connected by $4k-3$ internally vertex-disjoint paths, together
forming $T_0$.

Let $\WW=(W_1,\ldots,W_k)$ be a $k$-tuple of walks satisfying \rf{27no14b} and \rf{2fe15a}.
Consider a vertex $v$ of degree 2 in $T_0$, say incident with edges $e$ and $e'$ of $T_0$.
By \rf{2fe15a}, $W_i$ traverses edges in $\pi_e$ as often as it traverses edges
in $\pi_{e'}$.

For each $i=1,\ldots,k$, define $h_i:E(T_0)\to\{0,1,\ldots,2|E|\}$ by:
$h_i(e)$ is the number of times that $W_i$ traverses $\pi_e$, in any direction
(for $e\in E(T_0)$).
Then we can derive $\phi_{\WW}$ from $h_1,\ldots,h_k$, without knowing $\WW$:

\clnn{
Given $h_1,\ldots,h_k$, one can find $\phi_{\WW}$ in polynomial time.
}

\pfcl
Consider any $\{i,j\}\in F$.
Let $T''$ be the subgraph of $T'$ consisting of the edges traversed by $W_i$ and $W_j$.
We know $T''$ since we know $h_i$ and $h_j$.
We determine an undirected graph $H$ with vertex set $E(T'')$, calling two edges in $E(T'')$ {\em associated}
if they form an edge of $H$.

First, consider any vertex $v$ of $T_0$ of degree 2.
Let $e$ and $e'$ be the edges of $T_0$ incident with $v$, and consider the
parallel classes $\pi_e$ and $\pi_{e'}$ in $T''$.
As $|\pi_e|=h_i(e)+h_j(e)=h_i(e')+h_j(e')=|\pi_{e'}|$, we can order the edges in $\pi_e$ as
$e_1,\ldots,e_m$ from left to right when going towards $v$, and similarly, the edges in $\pi_{e'}$ as
$e'_1,\ldots,e'_m$, from left to right when going away from $v$.
For each $t=1,\ldots,m$, we `associate' $e_t$ and $e'_t$.

Next, consider any vertex $v$ of $T_0$ of degree 3.
Let $e$, $e'$, and $e''$ be the edges of $T_0$ incident with $v$.
Consider the parallel classes $\pi_e$, $\pi_{e'}$, and $\pi_{e''}$ in $T''$.
As $W_i$ and $W_j$ are non-returning (that is, satisfy \rf{2fe15a}(i)), we know that there exist nonnegative integers
$a$, $b$, and $c$ such that $|\pi_e|=b+c$, $|\pi_{e'}|=a+c$, and $|\pi_{e''}|=a+b$.
These numbers are unique and can be directly calculated from
$|\pi(e)|$, $|\pi(e')|$, and $|\pi(e'')|$.
This implies that the edges in $\pi_e\cup\pi_{e'}\cup\pi_{e''}$ can uniquely be pairwise `associated' in such a
fashion that any two associated pairs of edges are noncrossing at $v$ and
such that no two edges in the same parallel class are associated.

Finally, consider any vertex $v$ of $T_0$ of degree 1.
So $v$ belongs to $\{r_1,s_1,\ldots,r_k,s_k\}$.
Let $e$ be the edge of $T_0$ incident with $v$.
Let $e_1,\ldots,e_t$ be the edges in the parallel class $\pi_e$ of $T''$, in order.
`Associate' $e_i$ with $e_{t+1-i}$ for each $i=1,\ldots,\lfloor\frac12t\rfloor$.
So if $t$ is odd (which is the case if and only if $v\in\{r_i,s_i,r_j,s_j\}$),
one edge in $\pi_e$ remains unassociated at $v$ (namely the middle edge).

Then the graph $H$ with $E(T'')$ as vertex set and all pairs of associated edges of $T''$ as edges of $H$,
consists of two paths,
corresponding to $W_i$ and $W_j$ in $T'$.
These sets of edges form two walks that we can orient, one from $r_i$ to $s_i$, the other from
$r_j$ to $s_j$.
Then for each edge $e$ of $T_0$ we know the order, from left to right, in which $W_i$ and $W_j$ traverse the
parallel class $\pi_e$ of $T''$,
and we can derive the direction.
Concluding, we can derive the subword of $\phi_{\WW}(e)$ made up by the symbols $g_i$, $g_i^{-1}$,
$g_j$, and $g_j^{-1}$.
(It is important here that we know that $H$ comes from an $r_i-s_i$ walk $W_i$ and an
$r_j-s_j$ walk $W_j$.
So $H$ contains no circuit, for which we would not know whether it belongs to $W_i$ or to $W_j$.)

As we can do this for each $\{i,j\}\in F$, we can derive $\phi_{\WW}(e)$.
This follows from the fact that for any word $w$ with symbols $g_1,g_1^{-1},\ldots,g_k,g_k^{-1}$,
if we know for each $\{i,j\}\in F$ the subword $w_{i,j}$
of $w$ made up by $g_i$, $g_i^{-1}$, $g_j$, $g_j^{-1}$, we can determine $w$ as word up to transposition of
commuting symbols (but without cancellation):
Start by finding an $i\in[k]$ and $\alpha\in\{g_i,g_i^{-1}\}$ that occurs first in $w_{i,j}$ for each
$j$ with $\{i,j\}\in F$.
By transposition we can assume that $\alpha$ is the first symbol of $w$.
Then delete the first $\alpha$ from each $w_{i,j}$ with $\{i,j\}\in F$, and iterate.

Thus, temporarily, we do not cancel $g_i$ with $g_i^{-1}$ or $g_j$ with $g_j^{-1}$,
but work in the semigroup generated by $g_1,g_1^{-1},\ldots,g_k,g_k^{-1}$ with
relations $g_ig_j=g_jg_i$, $g_ig_j^{-1}=g_j^{-1}g_i$, and $g_i^{-1}g_j^{-1}=g_j^{-1}g_i^{-1}$ for
all distinct $i,j$ with $\{i,j\}\not\in F$.
At the end we factor out to the group $G_F$.

Concluding, we can find the element of $G_F$ represented by word $w$.
(Here we use the assumption that each $i$ is contained in some pair in $F$.)
\openbx

We finally describe the required algorithm.
Enumerate all $k$-tuples of functions $h_1,\ldots,h_k:E(T_0)\to\{0,1,\ldots,2|E|\}$ with the
property that if $e$ and $e'$ are the edges of $T_0$ incident with a vertex of $T_0$ of degree 2,
then $h_i(e)=h_i(e')$ for each $i$.
Determine, if possible, $\phi_{\WW}$.
All such $\phi_{\WW}$ form $\Phi$.

Since $T_0$ consists of vertices of degree 1 and 3 together with $4k-3$ 
internally vertex-disjoint paths connecting these vertices, there are $((2|E|+1)^{4k-3})^k$
such $k$-tuples $h_1,\dots,h_k$.
For fixed $k$, this is polynomially bounded.
\bx

\subsectz{Finding partially disjoint paths}

Concluding, we have:

\thmz{
For each fixed $k$, the partially disjoint paths problem in directed planar graphs is solvable
in polynomial time.
}

\pf
Directly from Propositions \ref{28no14b} and \ref{28no14c}.
\bx

\sectz{Some extensions and open questions}

The theorem can be extended to the case where for each edge $e$ of $D$ a subset $K_e$ of
$[k]$ is given, prescribing that $e$ may be traversed only by paths $P_i$ with $i\in K_e$.
This amounts to restricting $I$ in \rf{28no14d} to subsets of $K_e$.
Instead of requiring disjointness of certain pairs of paths, one may relax this to requiring that
certain pairs of paths are noncrossing: so they are allowed to `touch' each other in a vertex,
but not cross.
This amounts to deleting the `nonplanar' edges $e_{v,F,F'}$.

One may impose further conditions of the following kind.
Choose an (undirected) path $Q$ in the dual graph $D^*$, connecting two faces $F$ and $F'$ of $D$.
Then one may restrict the total `flow' of paths $P_i$ in $D$ that intersect $Q$: as long as the
restriction can be described by a closed subset of $G_F$, this requirement translates into an
extra nonplanar edge added to the dual graph $D^*$, like before we did for paths in $D^*$ connecting two
faces incident with a vertex $v$.

Moreover, the theorem extends to directed graphs $D$ on any fixed compact surface instead of
planar graphs.
Then, instead of considering the spanning tree $T$ in Section \ref{2de14c}.\ref{28no14e}, one considers a minimal
connected spanning subgraph $T$ that is cellularly embedded, i.e., each face is a disk
(assuming without loss of generality that $D$ is cellularly embedded).
Fixing the surface, the number of edges in $T$ is only a fixed amount more than in a spanning tree,
and the enumeration of homology classes can be bounded accordingly.

We finally mention some open questions.
The running time of the algorithm above is bounded by a polynomial with exponent depending on $k$
(in fact, $O(k^2)$).
This raises the question if the problem is `fixed parameter tractable'; that is, can the degree
of the polynomial be fixed independently of $k$, while the dependence of $k$ is only in the coefficient.
As mentioned, this question was answered confirmatively by Cygan, Marx, Pilipczuk, and Pilipczuk [3] for the $k$ {\em fully}
disjoint paths problem in directed planar graphs.

Another open question is if the condition of fixing $k$ can be relaxed
to fixing other parameters of the graph $\Gamma=([k],F)$.
One may think of fixing the maximum degree of $\Gamma$, or (more weakly) fixing the chromatic number
of $\Gamma$, or (even more weakly) fixing the clique number of $\Gamma$.
A different open question is if instead of fixing $k$, it suffices to fix the number of faces that
can cover all terminals (by the face boundaries).

A closely related open question is the complexity of the edge-disjoint version of the problem.
Even for the following problem it is not known whether it is polynomial-time solvable or NP-complete:
given a directed planar graph $D=(V,E)$ and vertices $r$ and $s$, find a directed $r-s$ path $P$
and a directed $s-r$ path $Q$ such that $P$ and $Q$ are edge-disjoint.
The corresponding problem for the undirected case is polynomial-time solvable (for any fixed number
of paths), even for general nonplanar graphs, by Robertson and Seymour [14].

Let us finally question whether the polynomial-time solvability of the cohomology feasibility problem
for graph groups (polynomial-time even for {\em unfixed} $k$) has other applications,
for instance to free partially commutative semigroups as studied for inhomogeneous sorting and scheduling
of concurrent processes (cf.\ Anisimov and Knuth [1], Diekert [4,\linebreak[0]5]).

\vspace{4mm}

\section*{References}\label{REF}
{\small
\begin{itemize}{}{
\setlength{\labelwidth}{4mm}
\setlength{\parsep}{0mm}
\setlength{\itemsep}{0mm}
\setlength{\leftmargin}{5mm}
\setlength{\labelsep}{1mm}
}
\item[\mbox{\rm[1]}] A.V. Anisimov, D.E. Knuth, 
Inhomogeneous sorting,
{\em International Journal of Computer and Information Sciences} 8 (1979) 255--260.

\item[\mbox{\rm[2]}] A. Baudisch, 
Kommutationsgleichungen in semifreien Gruppen,
{\em Acta Mathematica A\-ca\-de\-mi\-ae Scientiarum Hungaricae}
29 (1977) 235--249.

\item[\mbox{\rm[3]}] M. Cygan, D. Marx, M. Pilipczuk, M. Pilipczuk, 
The planar directed k-vertex-disjoint paths problem is fixed-parameter tractable,
in: {\em 2013 {IEEE} 54th Annual Symposium on Foundations of Computer Science ({FOCS})},
{IEEE}, 2013, pp. 197--206.

\item[\mbox{\rm[4]}] V. Diekert, 
On the Knuth-Bendix completion for concurrent processes,
{\em Theoretical Computer Science}  66 (1989) 117--136.

\item[\mbox{\rm[5]}] V. Diekert, 
{\em Combinatorics on Traces},
Lecture Notes in Computer Science 454,
Springer, Berlin, 1990.

\item[\mbox{\rm[6]}] C. Droms, 
Isomorphisms of graph groups,
{\em Proceedings of the American Mathematical Society} 100 (1987) 407--408.

\item[\mbox{\rm[7]}] E.S. Esyp, I.V. Kazachkov, V.N. Remeslennikov, 
Divisibility theory and complexity of algorithms for free partially commutative groups,
in: {\em Groups, Languages, Algorithms},
{\em Contemporary Mathematics} 378,
American Mathematical Society, Providence, R.I., 2005, pp. 319--348.

\item[\mbox{\rm[8]}] S. Fortune, J. Hopcroft, J. Wyllie, 
The directed subgraph homeomorphism problem,
{\em Theoretical Computer Science} 10 (1980) 111--121.

\item[\mbox{\rm[9]}] K. Kawarabayashi, Y. Kobayashi, B. Reed, 
The disjoint paths problem in quad\-rat\-ic time,
{\em Journal of Combinatorial Theory, Series B} 102 (2012) 424--435. 

\item[\mbox{\rm[10]}] J.F. Lynch, 
The equivalence of theorem proving and the interconnection problem,
{\em ({ACM}) {SIGDA} Newsletter} 5:3 (1975) 31--36.

\item[\mbox{\rm[11]}] R.C. Lyndon, P.E. Schupp, 
{\em Combinatorial Group Theory},
Springer, Berlin, 1977.

\item[\mbox{\rm[12]}] W. Magnus, A. Karrass, D. Solitar, 
{\em Combinatorial Group Theory},
Wiley-Inter\-science, New York, 1966.

\item[\mbox{\rm[13]}] B.A. Reed, N. Robertson, A. Schrijver, P.D. Seymour, 
Finding disjoint trees in planar graphs in linear time,
in: {\em Graph Structure Theory}
(N. Robertson, P. Seymour, eds.),
{\em Contemporary Mathematics} 147,
American Mathematical Society, Providence, R.I., 1993, pp. 295--301.

\item[\mbox{\rm[14]}] N. Robertson, P.D. Seymour, 
Graph minors. {XIII}. The disjoint paths problem,
{\em Journal of Combinatorial Theory, Series B} 63 (1995) 65--110.

\item[\mbox{\rm[15]}] A. Schrijver, 
Finding $k$ disjoint paths in a directed planar graph,
{\em {SIAM} Journal on Computing} 23 (1994) 780--788.

\item[\mbox{\rm[16]}] A. Schrijver, 
{\em Combinatorial Optimization --- Polyhedra and Efficiency},
Springer, Ber\-lin, 2003.

\item[\mbox{\rm[17]}] H. Servatius, 
Automorphisms of graph groups,
{\em Journal of Algebra} 126 (1989) 34--60.

\item[\mbox{\rm[18]}] C. Wrathall, 
The word problem for free partially commutative groups,
{\em Journal of Symbolic Computation} 6 (1988) 99--104.

\end{itemize}
}

\end{document}